\documentclass{amsart}
\usepackage{amssymb}

\newtheorem{theorem}{Theorem}[section]
\theoremstyle{plain}

\newtheorem{condition}{Condition}[section]

\newtheorem{corollary}{Corollary}[section]

\newtheorem{lemma}{Lemma}[section]

\newtheorem{problem}{Problem}[section]
\newtheorem{proposition}{Proposition}[section]
\newtheorem{remark}{Remark}[section]

\numberwithin{equation}{section}

\begin{document}
\title[\textbf{Fixed point indices and periodic points}]{\textbf{Fixed point
indices and periodic points of holomorphic mappings}}
\author{Guang yuan\medskip\ Zhang }
\address{Department of Mathematical Sciences, Tsinghua University, Beijing
100084, P. R. China}
\thanks{Project 10271063 and 10571009 supported by NSFC}
\email{gyzhang@math.tsinghua.edu.cn; gyzhang@mail.tsinghua.edu.cn}

\begin{abstract}
Let $\Delta ^{n}$ be the ball $|x|<1$ in the complex vector space $\mathbb{C}%
^{n}$, let $f:\Delta ^{n}\rightarrow \mathbb{C}^{n}$ be a holomorphic
mapping and let $M$ be a positive integer. Assume that the origin $%
0=(0,\dots ,0)$ is an isolated fixed point of both $f$ and the $M$-th
iteration $f^{M}$ of $f$. Then for each factor $m$ of $M,$ the origin is
again an isolated fixed point of $f^{m}$ and the fixed point index $\mu
_{f^{m}}(0)$ of $f^{m}$ at the origin is well defined, and so is the (local)
Dold's index (see [\ref{Do}]) at the origin:%
\begin{equation*}
P_{M}(f,0)=\sum_{\tau \subset P(M)}(-1)^{\#\tau }\mu _{f^{M:\tau }}(0),
\end{equation*}%
where $P(M)$ is the set of all primes dividing $M,$ the sum extends over all
subsets $\tau $ of $P(M)$, $\#\tau \ $is the cardinal number of $\tau $ and $%
M:\tau =M(\prod_{p\in \tau }p)^{-1}$.

$P_{M}(f,0)$ can be interpreted to be the number of periodic points of
period $M$ of $f$ overlapped at the origin: any holomorphic mapping $%
f_{1}:\Delta ^{n}\rightarrow \mathbb{C}^{n}$ sufficiently close to $f$ has
exactly $P_{M}(f,0)$ distinct periodic points of period $M$ near the origin$%
, $ provided that all the fixed points of $f_{1}^{M}$ near the origin are
simple. Note that $f$ itself has no periodic point of period $M$ near the
origin$.$

According to M. Shub and D. Sullivan's work [\ref{SS}], a necessary
condition so that $P_{M}(f,0)\neq 0$ is that the linear part of $f$ at the
origin has a periodic point of period $M.$ The goal of this paper is to
prove that this condition is sufficient as well for holomorphic mappings.
\end{abstract}

\subjclass[2000]{32H50, 37C25}
\maketitle

\section{\textbf{Introduction\label{S1}}}

\subsection{\textbf{Counting periodic points by the M\"{o}bius Inversion
Formula}\label{S1-1}}

To introduce Dold's index, we first present a formula for counting periodic
points. Let $E$ be a finite set and let $h$ be a mapping from $E$ into
itself. Then, for each positive integer $m$, the $m$-th iteration $h^{m}$ of
$h$ is well defined by $h^{0}=id,$ $h^{1}=h,\dots ,h^{m}=h\circ h^{m-1}$,
successively. An element $p\in E$ is called a \emph{periodic point} of $h$
with \emph{period} $m_{0}$ if $h^{m_{0}}(p)=p$ but $h^{m}(p)\neq p$ for $%
m=1,2,\dots ,m_{0}-1.$ Thus, the period is the smallest positive integer $%
m_{0}$ such that $h^{m_{0}}(p)=p.$

For each positive integer $m,$ we denote by $\mathrm{Fix}(h^{m})$ the set of
all fixed points of $h^{m}$, by $\mathcal{L}(h^{m})$ the cardinality of $%
\mathrm{Fix}(h^{m})$ and by $\mathcal{P}_{m}(h)$ the number of all periodic
points of $h$ with period $m.$ For a given positive integer $M,$ let us
derive the expression of the number $\mathcal{P}_{M}(h)$ in terms of the
numbers $\mathcal{L}(h^{m}).$

For each factor $M^{\prime }$ of $M$ and each $p\in \mathrm{Fix}%
(h^{M^{\prime }}),$ $p$ must be a periodic point of period $m$ for some
factor $m$ of $M^{\prime }$ (see Lemma \ref{ad-lem4-1})$.$ Thus%
\begin{equation}
\mathcal{L}(h^{M^{\prime }})=\sum_{m|M^{\prime }}\mathcal{P}_{m}(h),\
\mathrm{for\ all\ }M^{\prime }\in \mathbb{N}\ \mathrm{with}\text{\textrm{\ }}%
M^{\prime }|M,  \label{intro00}
\end{equation}%
where $\mathbb{N}$ is the set of all positive integers, the notation $%
M^{\prime }|M$ indicates that $M^{\prime }$ divides $M$ and, for each $%
M^{\prime },$ the sum extends over all factors $m$ of $M^{\prime }.$ Then we
can solve the system (\ref{intro00}) of linear equations with unknowns $%
\mathcal{P}_{m}(h),m|M.$ By the famous M\"{o}bius Inversion Formula (see
\cite{Hy}), the solution for $\mathcal{P}_{M}(h)$ can be expressed in terms
of the numbers $\mathcal{L}(h^{M^{\prime }})$ with $M^{\prime }|M$ as%
\begin{equation}
\mathcal{P}_{M}(h)=\sum_{\tau \subset P(M)}(-1)^{\#\tau }\mathcal{L}%
(h^{M:\tau }),  \label{intro0}
\end{equation}%
where $P(M)$ is the set of all primes dividing $M,$ the sum extends over all
subsets $\tau $ of $P(M),$ $\#\tau \ $is the cardinal number of $\tau $ and $%
M:\tau =M(\prod_{p\in \tau }p)^{-1}$. Note that the sum includes the term $%
\mathcal{L}(h^{M})$ which corresponds to the empty subset $\tau =\emptyset $%
. If $M=12=2^{2}\cdot 3,$ for example, then $P(M)=\{2,3\},$ and
\begin{equation*}
\mathcal{P}_{12}(h)=\mathcal{L}(h^{12})-\mathcal{L}(h^{4})-\mathcal{L}%
(h^{6})+\mathcal{L}(h^{2}).
\end{equation*}

\subsection{\textbf{Fixed point indices of holomorphic mappings}\label{S1-2}}

Let $\mathbb{C}^{n}$ be the complex vector space of dimension $n$, let $U$
be an open subset of $\mathbb{C}^{n}$ and let $g\in \mathcal{O}(U,\mathbb{C}%
^{n}),$ the space of holomorphic mappings from $U$ into $\mathbb{C}^{n}.$

If $p\in U$ is an isolated zero of $g,$ say, there exists a ball $B$
centered at $p$ with $\overline{B}\subset U$ such that $p$ is the unique
solution of the equation $g(x)=0$ ($0$ is the origin) in $\overline{B},$
then we can define the \emph{zero order} of $g$ at $p$ by%
\begin{equation*}
\pi _{g}(p)=\#(g^{-1}(q)\cap B)=\#\{x\in B;g(x)=q\},
\end{equation*}%
where $q$ is a regular value of $g$ such that $|q|$ is small enough and $\#$
denotes the cardinality. $\pi _{g}(p)$ is a well defined integer (see \cite%
{LL} or [\ref{Zh}] for the details).

If $q\ $is an isolated fixed point of $g,$ then $q$ is an isolated zero of
the mapping $id-g:U\rightarrow \mathbb{C}^{n},$ which puts each $x\in U$
into $x-g(x)\in \mathbb{C}^{n},$ and then the \emph{fixed point index} $\mu
_{g}(q)$ of $g$ at $q$ is defined to be the zero order of $id-g$ at $q:$%
\begin{equation*}
\mu _{g}(q)=\pi _{id-g}(q).
\end{equation*}

If $q$ is a fixed point of $g$ such that $id-g$ is regular at $q,$ say, the
Jacobian matrix $Dg(q)$ of $g$ at $q$ has no eigenvalue $1,$ $q$ is called a
\emph{simple} fixed point of $g.$ By Lemma \ref{cc3.2}, a fixed point of a
holomorphic mapping has index $1$ if and only if it is simple.

Actually, the zero order defined here is the (local) topological degree, and
the fixed point index defined here is the (local) Lefschetz fixed point
index, if $g$ is regarded as a continuous mapping of real variables. See the
Appendix section for the details.

\subsection{\textbf{The local Dold's indices for holomorphic mappings}}

Let $\Delta ^{n}$ be the ball $|x|<1$ in $\mathbb{C}^{n}$ and let $f\in
\mathcal{O}(\Delta ^{n},\mathbb{C}^{n})$. If the origin $0=(0,\dots ,0)$ is
a fixed point of $f,$ then for any positive integer $m,$ the $m$-th
iteration $f^{m}$ of $f$ is well defined in a neighborhood $V_{m}$ of the
origin$.$

Now, let $M$ be a positive integer and assume that the origin is an isolated
fixed point of both $f$ and $f^{M}$. Then for each factor $m$ of $M,$ the
origin is an isolated fixed point of $f^{m}$ as well and the fixed point
index $\mu _{f^{m}}(0)$ of $f^{m}$ at the origin is well defined. Therefore,
we can define the (local) \emph{Dold's index} [\ref{Do}] similar to (\ref%
{intro0}):%
\begin{equation}
P_{M}(f,0)=\sum_{\tau \subset P(M)}(-1)^{\#\tau }\mu _{f^{M:\tau }}(0).
\label{intro-0}
\end{equation}

The importance of this number is that it can be interpreted to be the number
of periodic points of $f$ of period $M$ overlapped at the origin$:$

\emph{For any ball }$B$\emph{\ centered at the origin, with }$\overline{B}%
\subset \Delta ^{n},$\emph{\ such that }$f^{M}$ \emph{has no fixed point in}
$\overline{B}\emph{\ other}$\emph{\ than the origin,} \emph{any }$f_{1}\in
\mathcal{O}(\Delta ^{n},\mathbb{C}^{n})$ \emph{has exactly} $P_{M}(f,0)$
\emph{mutually distinct periodic points of period} $M$ \emph{in} $B,$ \emph{%
provided that all fixed points of} $f_{1}^{M}$ \emph{in} $B$ \emph{are
simple and that} $f_{1}$ \emph{is sufficiently close to} $f,$ \emph{in the
sense that }%
\begin{equation*}
||f-f_{1}||_{\Delta ^{n}}=\sup_{x\in \Delta ^{n}}|f(x)-f_{1}(x)|
\end{equation*}
\emph{is small enough (see Corollary \ref{cc3.2+1} (iii) and Lemma \ref%
{Th3.0} (ii)). }

Note that $f$ itself has no periodic point of period $M$ in $\overline{B}.$
This gives rise to an interesting problem:

\begin{problem}
\label{Prob1}What is the condition under which $P_{M}(f,0)\neq 0?$
\end{problem}

By Corollary \ref{cc3.5+2}, this is equivalent to ask

\begin{problem}
\label{Prob2}What is the condition under which there exists a sequence of
holomorphic mappings $f_{j}:\Delta ^{n}\rightarrow \mathbb{C}^{n}$ such that

(1) $f_{j}$ uniformly converges to $f$ in a neighborhood of the origin$,$ and

(2) each $f_{j}$ has a periodic point of period $M$ converging to the origin
as $j\rightarrow \infty ?$
\end{problem}

According to M. Shub and D. Sullivan's work [\ref{SS}], a necessary
condition such that $P_{M}(f,0)\neq 0$ is that the linear part of $f$ at the
origin has a periodic point of period $M$ (see Lemma \ref{Th3.1} and its
consequence Lemma \ref{Th3+1} (i)). The term \textquotedblleft linear
part\textquotedblright\ indicates the linear mapping $l:\mathbb{C}%
^{n}\rightarrow \mathbb{C}^{n},$
\begin{equation*}
l(x_{1},\dots ,x_{n})=(\sum_{j=1}^{n}a_{1j}x_{j},\dots
,\sum_{j=1}^{n}a_{nj}x_{j}),
\end{equation*}%
where $(a_{ij})=Df(0)=(\frac{\partial f_{i}}{\partial x_{j}})|_{0}$ is the
Jacobian matrix of $f$ at the origin$.$

The goal of this paper is to prove that this condition is sufficient as well:

\begin{theorem}
\label{Th1.0}Let $M$ be a positive integer and let $f:\Delta ^{n}\rightarrow
\mathbb{C}^{n}$ be a holomorphic mapping such that the origin is an isolated
fixed point of both $f$ and $f^{M}$. Then
\begin{equation}
P_{M}(f,0)>0  \label{0000}
\end{equation}%
if and only if the linear part of $f$ at the origin has a periodic point of
period $M.$
\end{theorem}

When $M=1,$ this theorem is trivial, since any linear mapping has a fixed
point at the origin and by the assumption on $f$ and Lemma \ref{cc3.2}, $%
P_{1}(f,0)=\mu _{f}(0)>0$. When $M>1,$ by Lemma \ref{linear}, the conclusion
in Theorem \ref{Th1.0} can be restated as:

\begin{corollary}
$P_{M}(f,0)>0$ if and only if the following condition holds.
\end{corollary}

\begin{condition}
\label{1.0} There exist positive integers $m_{1},\dots ,m_{s}$ ($s\leq n$)
such that their least common multiple is $M$ and that $Df(0)$ has
eigenvalues that are primitive $m_{1}$-th, $\dots $, $m_{s}$-th roots of
unity, respectively.
\end{condition}

A complex number $\lambda $ is called a primitive $m$-th root of unity if $%
\lambda ^{m}=1$ but $\lambda ^{j}\neq 1$ for $j=1,2,\dots ,m-1.$ Simple
examples show that the sufficiency in Theorem \ref{Th1.0} fails when the
mapping $f$ is not holomorphic.

For positive integers $n_{1},\dots ,n_{k},$ their least common multiple will
be denoted by $[n_{1},\dots ,n_{k}].$ A linear mapping that has periodic
points of periods $n_{1},\dots ,n_{k}$ must have periodic points of period $%
[n_{1},\dots ,n_{k}].$ From this fact and Theorem \ref{Th1.0}, we can
conclude directly that if a holomorphic mapping has periodic points of
periods $n_{1},\dots ,n_{k}$ overlapped at a fixed point, then it has
periodic points of period $[n_{1},\dots ,n_{k}]$ overlapped at that fixed
point, say precisely, we have:

\begin{corollary}
\label{R1}Let $f:\Delta ^{n}\rightarrow \mathbb{C}^{n}$ be a holomorphic
mapping such that the origin is an isolated fixed point of $%
f,f^{n_{1}},\dots ,f^{n_{k}},$ where $n_{1},\dots ,n_{k}$ are positive
integers. If $P_{n_{j}}(f,0)>0$ for all $j=1,\dots ,k,$ and if the origin is
an isolated fixed point of $f^{[n_{1},\dots ,n_{k}]},$ then
\begin{equation*}
P_{[n_{1},\dots ,n_{k}]}(f,0)>0.
\end{equation*}
\end{corollary}

Now, let $M\in \mathbb{N}\backslash \{1\}$ ($\mathbb{N}$ is the set of all
positive integers) and let $f:\Delta ^{n}\rightarrow \mathbb{C}^{n}$ be a
holomorphic mapping such that the origin is an isolated fixed point of both $%
f$ and $f^{M},$ and that $Df(0)$ satisfies Condition \ref{1.0}. We shall
make some remarks on the inequality (\ref{0000}).

When $n=1,$ (\ref{0000}) is known, and can be deduced as follows. In this
case, Condition \ref{1.0} means that $Df(0)$ is a primitive $M$-th root of
unity, and then, for each $j=1,\dots ,M-1,$ the origin is a simple zero of
the one variable holomorphic function $x-f^{j}(x),$ which implies $\mu
_{f^{j}}(0)=1.$ Therefore, (\ref{intro-0}) becomes
\begin{eqnarray*}
P_{M}(f,0) &=&\mu _{f^{M}}(0)+\sum_{\substack{ \tau \subset P(M)  \\ \tau
\neq \emptyset }}(-1)^{\#\tau } \\
&=&\mu _{f^{M}}(0)+\sum_{k=1}^{\#P(M)}(-1)^{k}\left(
\begin{array}{c}
k \\
\#P(M)%
\end{array}%
\right) =\mu _{f^{M}}(0)-1,
\end{eqnarray*}%
and then (\ref{0000}) is equivalent to%
\begin{equation}
\mu _{f^{M}}(0)-1>0.  \label{01}
\end{equation}%
On the other hand, since $(Df^{M})(0)=\left( Df(0)\right) ^{M}=1,$ the
origin is a zero of the holomorphic function $x-f^{M}(x)$ with order at
least $2.$ Thus, (\ref{01}) holds. In fact, in this case $\mu
_{f^{M}}(0)=kM+1$ for some positive integer $k$ (see [\ref{Mil}]).

When $n=2,$ (\ref{0000}) follows from a result obtained by the author in [%
\ref{Zh}]: the sequence in Problem \ref{Prob2} exists when $Df(0)$ satisfies
Condition \ref{1.0}. For $n=2,$ Condition \ref{1.0} implies that either (i) $%
Df(0)$ has an eigenvalue $\lambda $ that is a primitive $M$-th root of
unity, or (ii) the two eigenvalues of $Df(0)$ are primitive $m_{1}$-th and $%
m_{2}$-th roots of unity, respectively$,$ and $M=[m_{1},m_{2}]>\max
\{m_{1},m_{2}\}$.

The inequality (\ref{0000}) is easy to prove in the first case (i): it can
be dealt with as the above one dimensional case by a small perturbation.

In the second case (ii), (\ref{0000}) is equivalent to%
\begin{equation}
\mu _{f^{M}}(0)-\mu _{f^{m_{1}}}(0)-\mu _{f^{m_{2}}}(0)+1>0.  \label{02}
\end{equation}%
In fact, in this case, $m_{1}\nmid m_{2},m_{2}\nmid m_{1},$ and then the
origin is a simple fixed point of $f$ with $\mu _{f}(0)=P_{1}(f,0)=1$ and by
Lemma \ref{linear}, each periodic point of the linear part of $f$ at the
origin only has period $1,$ $m_{1},m_{2}$ or $M=[m_{1},m_{2}].$ Therefore,
by Shub and Sullivan's work (see Lemmas \ref{Th3.1} and \ref{Th3+1} (i)), we
have%
\begin{equation*}
\mu _{f^{M}}(0)=P_{M}(f,0)+P_{m_{1}}(f,0)+P_{m_{2}}(f,0)+1,
\end{equation*}%
\begin{equation*}
\mu _{f^{m_{j}}}(0)=P_{m_{j}}(f,0)+1,j=1,2,
\end{equation*}%
and then we have
\begin{equation*}
P_{M}(f,0)=\mu _{f^{M}}(0)-\mu _{f^{m_{1}}}(0)-\mu _{f^{m_{2}}}(0)+1,
\end{equation*}%
and the equivalence of (\ref{0000}) and (\ref{02}) follows.

The inequality (\ref{02}) is obtained by the author in [\ref{Zh}], which is
the key ingredient in [\ref{Zh}] for solving Problem \ref{Prob2} for $n=2$.
When $n>2,$ the perturbation method used to prove inequalities such as (\ref%
{01}) and (\ref{02}), which strongly depends upon one variable complex
analysis, no longer works directly.

Fortunately, we can make up this shortage by involving an elementary result
of the normal form theory. In Section \ref{S3} we shall introduce some
results about fixed point indices of iterated holomorphic mappings for later
use. Most of them are known. In Section \ref{S4}, we shall compute Dold's
indices for mappings in a special case. Then, in Sections \ref{S5} and \ref%
{S6}, we shall show that the general case can be reduced to the special case
considered in Section \ref{S4}, by small perturbations and coordinate
transformations, and the inequality (\ref{0000}) will finally be proved in
the general case. The necessity in Theorem \ref{Th1.0} will be proved by the
way.

\section{\textbf{Some properties of fixed point indices of holomorphic
mappings}\label{S3}}

In this section we introduce some results about fixed point indices of
iterated holomorphic mappings for later use. Most of them are known.

Let $U$ be a bounded open subset of $\mathbb{C}^{n}$ and let $f\in \mathcal{%
O(}\overline{U},\mathbb{C}^{n}),$ the space of holomorphic mappings from $%
\overline{U}$ into $\mathbb{C}^{n}$. If $f$ has no fixed point on the
boundary $\partial U$, then the fixed point set \textrm{Fix}$(f)$ of $f$ is
a compact analytic subset of $U,$ and then it is finite (see [\ref{Ch}]);
and therefore, we can define the global fixed point index $L(f)$ of $f$ as:
\begin{equation*}
L(f)=\sum_{p\in \mathrm{Fix}(f)}\mu _{f}(p),
\end{equation*}%
which is just the number of all fixed points of $f$, counting indices. $L(f)$
is, in fact, the Lefschetz fixed point index of $f$ (see the Appendix
section for the details).

For each $m\in \mathbb{N},$ the $m$-th iteration $f^{m}$ of $f$ is
understood to be defined on
\begin{equation*}
K_{m}(f)=\cap _{k=0}^{m-1}f^{-k}(\overline{U})=\{x\in \overline{U}%
;f^{k}(x)\in \overline{U}\ \mathrm{for\ all\ }k=1,\dots ,m-1\},
\end{equation*}%
which is the largest set where $f^{m}$ is well defined. Since $U$ is
bounded, $K_{m}(f)$ is a compact subset of $\overline{U}$. Here, $f^{0}=id.$

Now, let us introduce the global Dold's index. Let $M\in \mathbb{N}$ and
assume that $f^{M}$ has no fixed point on the boundary $\partial U.$ Then,
for each factor $m$ of $M,$ $f^{m}$ again has no fixed point on $\partial U,$
and then the fixed point set $\mathrm{Fix}(f^{m})$ of $f^{m}$ is a compact
subset of $U,$ since $\mathrm{Fix}(f^{m})$ is a closed subset of $K_{m}(f)$
and $K_{m}(f)$ is a compact subset of $\overline{U}$. Thus, there exists an
open subset $V_{m}$ of $U$ such that $\mathrm{Fix}(f^{m})\subset
V_{m}\subset \overline{V_{m}}\subset U$ and $f^{m}$ is well defined on $%
\overline{V_{m}},$ and thus $L(f^{m}|_{\overline{V_{m}}})$ is well defined
and we write $L(f^{m})=L(f^{m}|_{\overline{V_{m}}}),$ where $f^{m}|_{%
\overline{V_{m}}}$ is the restriction of $f^{m}$ to $\overline{V_{m}}.$ In
this way, we can define the global Dold's index (see [\ref{Do}]) as (\ref%
{intro0}):%
\begin{equation}
P_{M}(f)=\sum_{\tau \subset P(M)}(-1)^{\#\tau }L(f^{M:\tau }).  \label{p0p}
\end{equation}

It is clear that, for any $m\in \mathbb{N}$ and any compact subset $K$ of $U$
with $\cup _{j=1}^{m}f^{j}(K)\subset U,$ there exists a neighborhood $%
V\subset U$ of $K,$ such that for any $g\in \mathcal{O(}\overline{U},\mathbb{%
C}^{n})$ sufficiently close to $f$ (in the sense that $||g-f||_{\overline{U}%
}=\max_{x\in \overline{U}}|g(x)-f(x)|$ is sufficiently small), the
iterations $g^{j},j=1,\dots ,m,$ are well defined on $\overline{V}$ and%
\begin{equation*}
||g-f||_{\overline{U}}\rightarrow 0\Longrightarrow \max_{1\leq j\leq
m}||g^{j}-f^{j}||_{\overline{V}}\rightarrow 0.
\end{equation*}%
We shall use these facts frequently and tacitly.\medskip

The following result follows from Rouch\'{e}'s theorem for equidimensional
holomorphic mappings directly.

\begin{lemma}[\protect\cite{LL}]
\label{cc3.1}Let $U$ be a bounded open subset of $\mathbb{C}^{n},$ let $f\in
\mathcal{O}(\overline{U},\mathbb{C}^{n})$ and assume that $f$ has no fixed
point on the boundary $\partial U.$ Then:

(1) $f$ has only finitely many fixed points in $U$ and for any $g\in
\mathcal{O}(\overline{U},\mathbb{C}^{n})$ that is sufficiently close to $f$ (%
$||f-g||_{\overline{U}}$ is small enough), $g$ has no fixed point on $%
\partial U,$ has only finitely many fixed points in $U$ and satisfies
\begin{equation*}
L(g)=\sum_{p\in \mathrm{Fix}(g)}\mu _{g}(p)=\sum_{p\in \mathrm{Fix}(f)}\mu
_{f}(p)=L(f);
\end{equation*}%
say, the number of fixed points of $g,$ counting indices, equals to that of $%
f$.

(2) In particular, if $p_{0}\in U$ is the unique fixed point of $f$ in $%
\overline{U}$ and if $||f-g||_{\overline{U}}$ is small enough, then
\begin{equation*}
L(g)=\sum_{p\in \mathrm{Fix}(g)}\mu _{g}(p)=\mu _{f}(p_{0})=L(f),
\end{equation*}%
and furthermore, if in addition all fixed points of $g$ are simple, then
\begin{equation*}
L(g)=\#\mathrm{Fix}(g)=\mu _{f}(p_{0}).
\end{equation*}
\end{lemma}

\begin{corollary}
\label{cc3.2+1}Let $M$ be a positive integer, let $U$ be a bounded open
subset of $\mathbb{C}^{n}$, let $f\in \mathcal{O}(\overline{U},\mathbb{C}%
^{n})$ and assume that $f^{M}$ has no fixed point on $\partial U$. Then:

(i). There exists an open subset $V$ of $U,$ such that $f^{M}$ is well
defined on $\overline{V},$ has no fixed point outside $V,$ and has only
finitely many fixed points in $V$.

(ii). For any $g\in \mathcal{O}(\overline{U},\mathbb{C}^{n})$ sufficiently
close to $f$, $g^{M}$ is well defined on $\overline{V},$ has no fixed point
outside $V$ and has only finitely many fixed points in $V;$ and furthermore,%
\begin{equation*}
L(g^{M})=L(f^{M})\mathrm{,}\text{ \ }P_{M}(g)=P_{M}(f).
\end{equation*}

(iii). In particular, if $p_{0}\in U$ is the unique fixed point of both $f$
and $f^{M}$ in $\overline{U},$ then for any $g\in \mathcal{O}(\overline{U},%
\mathbb{C}^{n})$ sufficiently close to $f$,
\begin{equation*}
L(g^{M})=L(f^{M})=\mu _{f^{M}}(p_{0})\mathrm{,}\text{ }%
P_{M}(g)=P_{M}(f)=P_{M}(f,p_{0}).
\end{equation*}
\end{corollary}

\begin{remark}
Under the assumption that $f^{M}$ has no fixed point on $\partial U,$ for
any factor $m$ of $M,$ the conclusions (i)--(iii) remain valid if $M$ is
replaced by $m,$ since $f^{m}$ has no fixed point on $\partial U$ as well.
\end{remark}

\begin{proof}
Recall that the domain of definition of $f^{M}$ is the set
\begin{equation*}
K_{M}(f)=\cap _{k=0}^{M-1}f^{-k}(\overline{U})=\{x\in \overline{U};\
f^{k}(x)\in \overline{U}\ \mathrm{for\ all\ }k=1,\dots ,M-1\},
\end{equation*}%
which is a compact subset of $\overline{U}$, where $f^{0}=id$. Then the
fixed point set $\mathrm{Fix}(f^{M})$ of $f^{M}$ equals to $\{x\in
K_{M}(f);f^{M}(x)=x\},$ which is a compact subset of $K_{M}(f).$ Therefore,
by the assumption that $f^{M}$ has no fixed point on $\partial U$, $\mathrm{%
Fix}(f^{M})$ is a compact subset of $U,$ and then there exists an open
subset $U^{\prime }$ of $U$ with $\overline{U^{\prime }}\subset U,$ such
that \textrm{Fix}$(f^{M})\subset U^{\prime }.$

Thus, for the open set $V=\cap _{k=0}^{M-1}f^{-k}(U^{\prime }),$ we conclude
that $f^{M}$ has no fixed point outside $V.$ On the other hand, since $%
\overline{V}\subset \cap _{k=0}^{M-1}f^{-k}(\overline{U^{\prime }})\subset
K_{M}(f),$ $f^{M}$ is well defined on $\overline{V}.$ Therefore, $f^{M}$ has
only finitely many fixed points in $V$ by Lemma \ref{cc3.1}, and hence (i)
holds for this open set $V$; and furthermore, we have
\begin{equation}
f^{k}(\overline{V})\subset \overline{U^{\prime }}\subset U,\mathrm{\ for\ all%
}\ k=1,\dots ,M-1.  \label{p7p}
\end{equation}

Now, let $g\in \mathcal{O}(\overline{U},\mathbb{C}^{n})$. To prove (ii), we
first prove:\medskip

(iv). If $||f-g||_{\overline{U}}$ is small enough, then $g^{M}$ is well
defined on $\overline{V}$ and has no fixed point outside $V.\medskip $

Clearly, $\overline{V}$ is compact. Therefore, by (\ref{p7p}), if $||f-g||_{%
\overline{U}}$ is small enough, then $g^{k}(\overline{V})\subset U$ for all $%
k=1,\dots ,M-1,$ and then $g^{M}$ is well defined on $\overline{V}.$

Thus, if (iv) fails, then there exists a sequence $\{g_{j}\}\subset \mathcal{%
O}(\overline{U},\mathbb{C}^{n}),$ uniformly converging to $f,$ such that for
each $j,$ $g_{j}^{M}$ is well defined on $\overline{V}$ and $g_{j}^{M}$ has
a fixed point $p_{j}\in \overline{U}\backslash V.$ Since $\overline{U}%
\backslash V$ is compact, we may assume that%
\begin{equation}
p_{j}\rightarrow p,\ \mathrm{as\ }j\rightarrow \infty ,  \label{q2q}
\end{equation}%
for some $p\in \overline{U}\backslash V.$

But then, we can prove a contradiction: $f^{M}(p)$ is well defined and%
\begin{equation}
f^{M}(p)=p,  \label{q3q}
\end{equation}%
which contradicts the proved conclusion that $f^{M}$ has no fixed point
outside $V.$

What $g_{j}^{M}$ has a fixed point $p_{j}$ means that $g_{j}^{M}$ is well
defined at $p_{j},$ say,
\begin{equation}
g_{j}^{k}(p_{j})\in \overline{U}\mathrm{\ for\ all\ }k=1,\dots ,M,
\label{p1p}
\end{equation}%
and
\begin{equation}
g_{j}^{M}(p_{j})=p_{j}.  \label{p2p}
\end{equation}

We first show that for each $k\in \mathbb{N}$ with $k\leq M$, $f^{k}(p)$ is
well defined and%
\begin{equation}
f^{k}(p)\in \overline{U}\mathrm{\ and\ }|f^{k}(p)-g_{j}^{k}(p_{j})|%
\rightarrow 0\mathrm{\ as\ }j\rightarrow \infty .  \label{p3p}
\end{equation}

We prove this by induction on $k$. By (\ref{q2q}), by the continuity of $f$
at $p,$ and by the uniform convergence of $g_{j}$ (note that $f(p)$ is well
defined), we have
\begin{equation*}
|f(p)-g_{j}(p_{j})|\leq |f(p)-f(p_{j})|+|f(p_{j})-g_{j}(p_{j})|\rightarrow
0,\ \mathrm{as\ }j\rightarrow \infty .
\end{equation*}%
Then by (\ref{p1p}) we have $f(p)\in \overline{U}.$ Therefore, the
conclusion is true for $k=1.$

Let $t$ be any positive integer with $t<M$ and assume that the conclusion is
valid for $k=t.$ Then $f^{t}(p)\in \overline{U},$ and then by (\ref{p1p}),
both $f^{t+1}(p)$ and $f(g_{j}^{t}(p_{j}))$ are well defined. Therefore,
again by the induction hypothesis, by the continuity of $f$ at $f^{t}(p)$
and by the uniform convergence of $g_{j}$ we have
\begin{eqnarray*}
|f^{t+1}(p)-g_{j}^{t+1}(p_{j})| &\leq
&|f(f^{t}(p))-f(g_{j}^{t}(p_{j}))|+|f(g_{j}^{t}(p_{j}))-g_{j}(g_{j}^{t}(p_{j}))|
\\
&\rightarrow &0\mathrm{\ as\ }j\rightarrow \infty ,
\end{eqnarray*}%
which implies $f^{t+1}(p)\in \overline{U}$ by (\ref{p1p})$.$ Therefore, the
conclusion is true for $k=t+1.$ This completes the induction.

By (\ref{p3p}), $f^{M}(p)$ is well defined and $|f^{M}(p)-g_{j}^{M}(p_{j})|$
tends to $0$ as $j$ tends to $\infty ,$ and then (\ref{q3q}) holds by (\ref%
{q2q}) and (\ref{p2p}). This completes the proof of (iv).

By (iv) and Lemma \ref{cc3.1}, if $||f-g||_{\overline{U}}$ is small enough,
then the set of all fixed points of $g^{M}$ is contained in $V$ and is
finite.

It is clear that (iv) remains valid if $M$ is replaced by any factor $m$ of $%
M.$ On the other hand, if $||f-g||_{\overline{U}}$ is small enough$,$ then $%
||g^{m}-f^{m}||_{\overline{V}}$ is also small enough for all factors $m$ of $%
M.$ Therefore, if $||f-g||_{\overline{U}}$ is small enough, then by (i),
(iv) and Lemma \ref{cc3.1}, we have
\begin{equation*}
L(f^{m})=L(f^{m}|_{\overline{V}})=L(g^{m}|_{\overline{V}})=L(g^{m}),\
\mathrm{for\ all\ factors\ }m\ \mathrm{of\ }M,
\end{equation*}%
and then $P_{M}(f)=P_{M}(g)$ by the formula (\ref{p0p}), and (ii) is proved.

By the hypothesis in (iii) and the definitions, we have $\mu
_{f^{M}}(p_{0})=L(f^{M})$ and $P_{M}(f,p_{0})=P_{M}(f),$ and then (iii)
follows from (ii) directly.
\end{proof}

\begin{lemma}[\protect\cite{LL}]
\label{cc3.2}Let $U$ be an open subset of $\mathbb{C}^{n}$, let $f\in
\mathcal{O}(U,\mathbb{C}^{n})$ and assume that $p\in U$ is an isolated fixed
point of $f.$ Then $\mu _{f}(p)\geq 1,$ and equality holds if and only if $p$
is a simple fixed point of $f$, say, $Df(p)$ has no eigenvalue $1$.
\end{lemma}

\begin{corollary}
\label{adc}Let $U$ be a bounded open subset of $\mathbb{C}^{n},$ let $f\in
\mathcal{O}(\overline{U},\mathbb{C}^{n})$ and assume that $f^{M}$ has no
fixed point on $\partial U.$ If $f$ has a periodic point $p\in U$ with
period $M$, then any $g\in \mathcal{O}(\overline{U},\mathbb{C}^{n})$ that is
sufficiently close to $f$ has a periodic point with period $M$ in $U.$
\end{corollary}

\begin{proof}
By the assumption, $p$, $f(p),\dots ,f^{M-1}(p)$ are distinct each other and
all located in $U,$ and moreover, by Corollary \ref{cc3.2+1} (i), $p$ is an
isolated fixed point of $f^{M}$. Therefore, there exists a ball $B$ in $U$
centered at $p$ such that $p$ is the unique fixed point of $f^{M}$ in $%
\overline{B}$ and, for any $g\in \mathcal{O}(\overline{U},\mathbb{C}^{n})$
such that $||g-f||_{\overline{U}}$ is small enough$,$ $g^{M}$ is well
defined on $\overline{B}$ and has no fixed point on $\partial B$, and
\begin{equation}
\overline{B}\cap g^{k}(\overline{B})=\emptyset ,1\leq k\leq M-1.
\label{adc1}
\end{equation}%
Therefore, by Lemmas \ref{cc3.1} and \ref{cc3.2}, if $||g-f||_{\overline{U}}$
is small enough, then
\begin{equation*}
\sum_{q\in \mathrm{Fix}(g^{M}|_{\overline{B}})}\mu _{g^{M}}(q)=L(g^{M}|_{%
\overline{B}})=L(f^{M}|_{\overline{B}})=\mu _{f^{M}}(p)\geq 1,
\end{equation*}%
and then $g^{M}$ has a fixed point $q$ in $B,$ and by (\ref{adc1}), $q$ is a
periodic point of $g$ with period $M.$
\end{proof}

A fixed point $p$ of $f$ is called \textit{hyperbolic} if $Df(p)$ has no
eigenvalue of absolute $1.$ If $p$ is a hyperbolic fixed point of $f,$ then
it is a hyperbolic fixed point of all iterations $f^{j},j\in \mathbb{N}.$

\begin{lemma}
\label{cc3.3}Let $M$ be a positive integer, let $U$ be a bounded open subset
of $\mathbb{C}^{n}$and let $f\in \mathcal{O}(\overline{U},\mathbb{C}^{n}).$
Assume that $V$ is an open subset of $U$ such that $f^{M}$ is well defined
on $\overline{V}$ and has no fixed point on $\partial V$. Then there exists
a sequence $\{f_{j}\}\subset \mathcal{O}(\overline{U},\mathbb{C}^{n}),$
uniformly converging to $f,$ such that for each $j,$ $f_{j}^{M}$ is well
defined on $\overline{V}$ and all the fixed points of $f_{j}^{M}$ located in
$\overline{V}$ are hyperbolic.
\end{lemma}

A proof of this result follows from the argument in [\ref{B}]. Another proof
can be found in [\ref{Zh1}].

\begin{corollary}
\label{ad-c-1}Let $M$ be a positive integer, let $U$ be a bounded open
subset of $\mathbb{C}^{n}$, let $f\in \mathcal{O}(\overline{U},\mathbb{C}%
^{n})$ and assume that $f^{M}$ has no fixed point on $\partial U$. Then
there exists a sequence $\{f_{j}\}\subset \mathcal{O}(\overline{U},\mathbb{C}%
^{n})$, uniformly converging to $f,$ such that $f_{j}^{M}$ has only finitely
many fixed points, all of which are hyperbolic, for each $j.$
\end{corollary}

\begin{proof}
By the hypothesis, Corollary \ref{cc3.2+1} applies. Thus, the conclusion
follows from applying Lemma \ref{cc3.3} to $f:\overline{U}\rightarrow
\mathbb{C}^{n}$ and the set $V$ given by Corollary \ref{cc3.2+1}.
\end{proof}

\begin{lemma}
\label{Th3.0}Let $M$ be a positive integer, let $U$ be a bounded open subset
of $\mathbb{C}^{n}$ and let $f\in \mathcal{O}(\overline{U},\mathbb{C}^{n}).$
Assume that $f^{M}$ has no fixed point on $\partial U$ and each fixed point
of $f^{M}$ is simple. Then $\mathrm{Fix}(f^{M})$ is finite, and

\textrm{(i). }$L(f^{M})=\#\mathrm{Fix}(f^{M})=\sum_{m|M}P_{m}(f).$

\textrm{(ii). }$P_{M}(f)$ is the cardinality of the set of periodic points
of $f$ with period $M$.
\end{lemma}

\begin{proof}
This was proved in [\ref{FL}]. In fact, each fixed point of $f^{m}$ with $%
m|M $ is also simple, and then it has index $1.$ On the other hand, by
Corollary \ref{cc3.2+1} (i), for each factor $m$ of $M,$ $\mathrm{Fix}%
(f^{m}) $ is finite, and then $L(f^{m})$ is equal to the cardinality of $%
\mathrm{Fix}(f^{m}).$ Therefore, (i) and (ii) follows from the argument in
Section \ref{S1-1} (see also [\ref{Do}], p. 421--422).
\end{proof}

\begin{corollary}
\label{cc3.5+2}Let $M$ be a positive integer and let $f\in \mathcal{O}(%
\overline{\Delta ^{n}},\mathbb{C}^{n}).$ Assume that the origin is an
isolated fixed point of both $f$ and $f^{M}.$ Then $P_{M}(f,0)\geq 0,$ and
furthermore, the following three conditions are equivalent.

(a). $P_{M}(f,0)>0.$

(b). There exist a ball $B\subset \Delta ^{n}$ centered at the origin and a
sequence $\{f_{j}\}\subset \mathcal{O}(\overline{B},\mathbb{C}^{n})$
uniformly converging to $f$ on $\overline{B},$ such that for each $f_{j},$
\begin{equation}
f_{j}(0)=0\;\mathrm{and\;}P_{M}(f_{j},0)>0.  \label{apr2-3}
\end{equation}

(c). There exist a ball $B\subset \Delta ^{n}$ centered at the origin and a
sequence $\{f_{j}\}\subset \mathcal{O}(\overline{B},\mathbb{C}^{n})$
uniformly converging to $f$ on $\overline{B},$ such that each $f_{j}$ has a
periodic point of period $M$ converging to the origin$.$
\end{corollary}

By the following proof, the ball $B$ in (b) and (c) can be replaced with $%
\Delta ^{n}.$

\begin{proof}
By the hypothesis, there exists a ball $B_{1}\subset \Delta ^{n}$ centered
at the origin such that $f^{M}$ is well defined on $\overline{B_{1}}$ and%
\begin{equation}
\overline{B_{1}}\cap \mathrm{Fix}(f)=\overline{B_{1}}\cap \mathrm{Fix}%
(f^{M})=\{0\}.  \label{apr4-1}
\end{equation}%
Then by the definition we have
\begin{equation}
P_{M}(f,0)=P_{M}(f|_{\overline{B_{1}}}),  \label{apr3-1}
\end{equation}%
and, by Corollary \ref{cc3.2+1} and Lemma \ref{cc3.3}, there exists a $g\in
\mathcal{O}(\overline{\Delta ^{n}},\mathbb{C}^{n}),$ which can be chosen
arbitrarily close to $f,$ such that $g^{M}$ is well defined on $\overline{%
B_{1}}$, all fixed points of $g^{M}$ located in $\overline{B_{1}}$ are
hyperbolic and%
\begin{equation}
\partial B_{1}\cap \mathrm{Fix}(g^{M})=\emptyset \mathrm{\ and\ }P_{M}(g|_{%
\overline{B_{1}}})=P_{M}(f|_{\overline{B_{1}}}).  \label{apr-1}
\end{equation}%
By Lemma \ref{Th3.0}, we have $P_{M}(g|_{\overline{B_{1}}})\geq 0$, and then
we have $P_{M}(f,0)\geq 0.$ It remains to prove the equivalences.

If $P_{M}(f,0)>0,$ then we have by (\ref{apr3-1}) and (\ref{apr-1}) that $%
P_{M}(g|_{\overline{B_{1}}})>0,$ and then again by Lemma \ref{Th3.0}, $g$
has a periodic point $x_{g}$ of period $M$ in $B_{1},$ which converges to
the origin as $g$ uniformly converges to $f,$ by (\ref{apr4-1}). Thus, (%
\emph{a}) implies (\emph{c}).

If there exist a ball $B\subset \Delta ^{n}$ and a sequence $%
\{f_{j}\}\subset \mathcal{O}(\overline{B},\mathbb{C}^{n})$ satisfying (\emph{%
b}) or (\emph{c}), then we may assume $B=B_{1}$, by shrinking $B$ or $B_{1}.$
Then by (\ref{apr4-1}), (\ref{apr3-1}) and Corollary \ref{cc3.2+1}, for a
sufficiently large $j=j_{0},$ we have
\begin{equation}
\partial B_{1}\cap \mathrm{Fix}(f_{j_{0}}^{M})=\emptyset \text{\textrm{\
and\ }}P_{M}(f,0)=P_{M}(f|_{\overline{B_{1}}})=P_{M}(f_{j_{0}}),
\label{apr3-2}
\end{equation}%
and then by Corollaries \ref{cc3.2+1} and \ref{ad-c-1}, there exists an $%
h\in \mathcal{O}(\overline{B_{1}},\mathbb{C}^{n})$, which can be chosen
arbitrarily close to $f_{j_{0}},$ such that
\begin{equation}
\partial B_{1}\cap \mathrm{Fix}(h^{M})=\emptyset \text{\textrm{\ and\ }}%
P_{M}(h)=P_{M}(f_{j_{0}}),  \label{apr3-3}
\end{equation}%
and that all fixed points of $h^{M}$ are hyperbolic.

If (\emph{c}) holds, then by Corollary \ref{adc}, $h$ may be taken so close
to $f_{j_{0}}$ that $h$ has a periodic point of period $M,$ and then by
Lemma \ref{Th3.0}, $P_{M}(h)>0.$ Therefore, by (\ref{apr3-2}) and (\ref%
{apr3-3}), we have $P_{M}(f,0)>0$, and hence, (\emph{c}) implies (\emph{a}).

By the first part of (\ref{apr3-2}) and Corollary \ref{cc3.2+1}$,$ there
exists a ball $B_{2}\subset B_{1}$ centered at the origin such that $%
f_{j_{0}}^{M}$ has no fixed point in $\overline{B_{2}}$ other than the
origin. Therefore, if (\emph{b}) holds, then again by Corollary \ref{cc3.2+1}%
, we may take $h$ so close to $f_{j_{0}}$ that
\begin{equation*}
P_{M}(h|_{\overline{B_{2}}})=P_{M}(f_{j_{0}}|_{\overline{B_{2}}%
})=P_{M}(f_{j_{0}},0)>0,
\end{equation*}%
and all the above properties of $h$ remain. But then $h$ has a periodic
point in $B_{2}$ of period $M,$ by Lemma \ref{Th3.0}. Since all fixed points
of $h^{M}$ are hyperbolic, again by Lemma \ref{Th3.0}, we have $P_{M}(h)>0,$
and then by (\ref{apr3-2}) and (\ref{apr3-3}) we have $P_{M}(f,0)>0,$ and
then (\emph{b}) implies (\emph{a}).

It is trivial that (\emph{a}) implies (\emph{b}). This completes the proof.
\end{proof}

\begin{lemma}
\label{linear}Let $L:\mathbb{C}^{n}\rightarrow \mathbb{C}^{n}$ be a linear
mapping and let $M>1\ $be a positive integer. Then $L$ has a periodic point
of period $M$ if and only if there exist positive integers $m_{1},\dots
,m_{s}$ ($s\leq n$) such that $M=[m_{1},\dots ,m_{s}]$ and that $L$ has
eigenvalues that are primitive $m_{1}$-th, $\dots $, $m_{s}$-th roots of
unity, respectively.
\end{lemma}

This is a basic knowledge of elementary linear algebra. Recall that $%
[m_{1},\dots ,m_{s}]$ denotes the least common multiple of $m_{1},\dots
,m_{s}$. The following result is due to M. Shub and D. Sullivan [\ref{SS}].
It is also proved in [\ref{Zh1}].

\begin{lemma}
\label{Th3.1}Let $m>1$ be a positive integer and let $f\in \mathcal{O}%
(\Delta ^{n},\mathbb{C}^{n}).$ Assume that the origin is an isolated fixed
point of $f$ and that, for each eigenvalue $\lambda $ of $Df(0),$ either $%
\lambda =1$ or $\lambda ^{m}\neq 1$. Then the origin is still an isolated
fixed point of $f^{m}$ and%
\begin{equation*}
\mu _{f}(0)=\mu _{f^{m}}(0).
\end{equation*}
\end{lemma}

\begin{corollary}
\label{cc3.6}Let $m>1$ be a positive integer and let $f\in \mathcal{O}(%
\overline{\Delta ^{n}},\mathbb{C}^{n}).$ Assume that the origin is an
isolated fixed point of both $f$ and $f^{m},$ and that the linear part of $f$
at the origin has no periodic point of period $m$. Then there exists a ball $%
B$ centered at the origin$,$ such that for any $g\in \mathcal{O}(\overline{%
\Delta ^{n}},\mathbb{C}^{n})$ that is sufficiently close to $f$ on $%
\overline{B},$ $g$ has no periodic point of period $m$ in $B$.
\end{corollary}

\begin{proof}
This was proved in [\ref{CMY}] and [\ref{Zh}] for $C^{1}$ mappings. Here we
give a much simpler proof. If the linear part of $f$ at the origin has no
periodic point of period $m,$ then by Lemma \ref{linear}, for the least
common multiple $m^{\ast }$ of all numbers of the set%
\begin{equation*}
\{m^{\prime }\in \mathbb{N};m^{\prime }|m\mathrm{\ and\ }Df(0)\mathrm{\ has\
an\ eigenvalue\ that\ is\ a\ primitive\ }m^{\prime }\text{-}\mathrm{th\
root\ of\ }1\},
\end{equation*}%
we have $m^{\ast }|m$ and $m^{\ast }<m$ (if the above set is empty, then
write $m^{\ast }=1$), and each eigenvalue $\lambda $ of $Df(0)$ with $%
\lambda ^{m}=1$ satisfies $\lambda ^{m^{\ast }}=1.$ Then for $d=m/m^{\ast }$
and for any eigenvalue $\eta $ of $Df^{m^{\ast }}(0),$ either $\eta =1$ or $%
\eta ^{d}\neq 1,$ which implies that for each factor $d_{1}$ of $d,$ either $%
\eta =1$ or $\eta ^{d_{1}}\neq 1.$ Therefore, by the above lemma, $\mu
_{f^{m^{\ast }}}(0)=\mu _{(f^{m\ast })^{d_{1}}}(0)$ for any factor $d_{1}$
of $d,$ and then%
\begin{equation*}
P_{d}(f^{m^{\ast }},0)=\sum_{\tau \subset P(d)}(-1)^{\#\tau }\mu _{(f^{m\ast
})^{d:\tau }}(0)=\mu _{f^{m\ast }}(0)\sum_{\tau \subset P(d)}(-1)^{\#\tau
}=0.
\end{equation*}

If the conclusion of the corollary fails, then there exists a sequence $%
\{f_{j}\}\subset \mathcal{O}(\overline{\Delta ^{n}},\mathbb{C}^{n})$,
uniformly converging to $f$ in a neighborhood of the origin, such that $%
f_{j} $ has a periodic point $x_{j}$ with period $m$ converging to the origin%
$,$ and then $f_{j}^{m^{\ast }}$ uniformly converges to $f^{m^{\ast }}$ in a
neighborhood of the origin and $x_{j}$ is a periodic point of $%
f_{j}^{m^{\ast }}$ with period $d$ converging to the origin, and then by
Corollary \ref{cc3.5+2} we have $P_{d}(f^{m^{\ast }},0)>0.$ This is a
contradiction, and the proof is complete.
\end{proof}

\begin{lemma}
\label{Th3+1}Let $f\in \mathcal{O}(\overline{\Delta ^{n}},\mathbb{C}^{n})$
with $f(0)=0,$\ and let
\begin{equation*}
\mathfrak{M}_{f}=\{m\in \mathbb{N};\;\mathrm{the\;linear\;part\;of\;}f\;%
\mathrm{at}\;0\;\mathrm{has\;a\ periodic\;point\;of\;period\;}m\}.
\end{equation*}%
Then:

\textrm{(i). }For each $m\in \mathbb{N}\backslash \mathfrak{M}_{f}$ such
that the origin is an isolated fixed point of $f^{m},$%
\begin{equation*}
P_{m}(f,0)=0.
\end{equation*}

\textrm{(ii). }For each positive integer $M$ such that the origin is an
isolated fixed point of $f^{M},$%
\begin{equation}
\mu _{f^{M}}(0)=\sum_{m\in \mathfrak{M}_{f},\ m|M}P_{m}(f,0).  \label{gre}
\end{equation}
\end{lemma}

\begin{proof}
This was essentially proved in [\ref{CMY}] in a more complicated setting for
$C^{1}$ mappings. Here we give a much simpler proof.

(i) follows from Corollaries \ref{cc3.5+2} and \ref{cc3.6} directly. To
prove (ii), let $M\;$be a positive integer such that the origin is an
isolated fixed point of $f^{M}.$ Since all quantities in (\ref{gre}) are
completely determined by the local behavior of $f$ near the origin$,$ by
shrinking $\Delta ^{n},$ we may assume that the origin is the unique fixed
point of both $f$ and $f^{M}.$ Then, the origin is the unique fixed point of
$f^{m}$ for each factor $m$ of $M,$ and then by Corollaries \ref{cc3.2+1}
(iii) and \ref{cc3.3}, there exists a holomorphic mapping $g:\overline{%
\Delta ^{n}}\rightarrow \mathbb{C}^{n}$ sufficiently close to $f$, such that
$g^{M}$ has no fixed point on $\partial \Delta ^{n}$, all fixed points of $%
g^{M}$ are hyperbolic and
\begin{equation}
\mu _{f^{M}}(0)=L(g^{M}),\ P_{m}(f,0)=P_{m}(g)\mathrm{\;}\text{\textrm{for }}%
\mathrm{each\;}m|M.  \label{h1}
\end{equation}

Now, Lemma \ref{Th3.0} applies to $g,$ and we have%
\begin{equation*}
L(g^{M})=\sum_{m|M}P_{m}(g).
\end{equation*}%
Therefore, by (\ref{h1}) and (i) we have%
\begin{equation*}
\mu _{f^{M}}(0)=\sum_{m|M}P_{m}(g)=\sum_{m|M}P_{m}(f,0)=\sum_{m\in \mathfrak{%
M}_{f},\ m|M}P_{m}(f,0).
\end{equation*}%
This completes the proof.
\end{proof}

It is well known that the fixed point index is invariant after any
biholomorphic transformation of coordinate, and therefore, so is the Dold's
index, by the definition. Thus we have:

\begin{lemma}
\label{cc3.4}Let $f,h\in \mathcal{O}(\Delta ^{n},\mathbb{C}^{n})$ and let $k$
be a positive integer such that the origin is an isolated fixed point of
both $f$ and $f^{k}$. Assume $h(0)=0$ and $\det Dh(0)\neq 0$ and write $%
g=h\circ f\circ h^{-1}.$ Then the origin is still an isolated fixed point of
both $g$ and $g^{k}$, $\mu _{f^{k}}(0)=\mu _{g^{k}}(0)$ and $%
P_{k}(f,0)=P_{k}(g,0).$
\end{lemma}

\section{\textbf{Computations of Dold's indices under Conditions} \protect
\ref{c1}--\protect\ref{c5}\label{S4}}

In this section, we shall apply Cronin's theorem to compute Dold's indices
(at the origin) of the holomorphic mapping $f=(f_{1},\dots ,f_{n}):\Delta
^{n}\rightarrow \mathbb{C}^{n}$ given by%
\begin{equation}
\left\{
\begin{array}{ll}
f_{j}(x_{1},\dots ,x_{n})=\smallskip & \lambda
_{j}x_{j}+x_{j}\sum_{i=1}^{s}a_{ji}x_{i}^{m_{i}}+o_{j},\;j=1,\dots
,s,\smallskip \\
f_{r}(x_{1},\dots ,x_{n})= & \mu _{r}x_{r}+R_{r},\;r=s+1,\dots ,n,%
\end{array}%
\right.  \label{form-1}
\end{equation}%
associated with the following five conditions.\medskip

\begin{condition}
\label{c1}$\lambda _{1},\dots ,\lambda _{s}$ are primitive $m_{1}$-th$,\dots
,m_{s}$-th roots of unity, respectively, and $m_{1},\dots ,m_{s}$ are
mutually distinct primes. Thus, $m_{1}\geq 2,\dots ,m_{s}\geq 2.$
\end{condition}

\begin{condition}
\label{c2}$\mu _{r}^{m_{1}\dots m_{s}}\neq 1$ for each $r=s+1,\dots ,n.$
\end{condition}

\begin{condition}
\label{c3}Each principal submatrix\footnote{%
For the $n\times n$ matrix $A=(a_{ij})$, a $k\times k$ submatrix $%
B=(b_{st})=(a_{i_{s}j_{t}})$ of $A$ is called a principal submatrix of $A$
if it is obtained from $A$ by deleting some $n-k$ rows and deleting the
columns in the same order, say, $i_{1}=j_{1},\dots ,i_{k}=j_{k}.$} of $%
A=(a_{ji})$ is invertible$.$
\end{condition}

\begin{condition}
\label{c4}Each $o_{j}$ is a power series of the form
\begin{equation*}
o_{j}(x_{1},\dots ,x_{n})=x_{j}h_{j}(x_{1}^{m_{1}},\dots
,x_{s}^{m_{s}})+o_{j}^{\prime }\left( x_{1},\dots ,x_{n}\right) ,
\end{equation*}%
such that $h_{j}(x_{1}^{m_{1}},\dots ,x_{s}^{m_{s}})$ is a polynomial in $%
x_{1}^{m_{1}},\dots ,x_{s}^{m_{s}}$ of the form%
\begin{equation}
h_{j}(x_{1}^{m_{1}},\dots ,x_{s}^{m_{s}})=\sum_{i_{1}+i_{2}+\dots
+i_{s}=2}^{N_{j}}c_{i_{1}i_{2}\dots
i_{s}}^{j}(x_{1}^{m_{1}})^{i_{1}}(x_{2}^{m_{2}})^{i_{2}}\dots
(x_{s}^{m_{s}})^{i_{s}},  \label{ad2}
\end{equation}%
in which $N_{j}$ is a positive integer and all $c_{i_{1}i_{2}\dots
i_{s}}^{j} $ are constants, and that $o_{j}^{\prime }\left( x_{1},\dots
,x_{n}\right) $ is a power series in $x_{1},\dots ,x_{n},$ consisting of
terms of degree larger than $\left( m_{1}\dots m_{s}\right) ^{2}$.
\end{condition}

\begin{condition}
\label{c5}Each $R_{r}$ is a power series in $x_{1},\dots ,x_{n}$ consisting
of terms of degree $\geq 2$.
\end{condition}

Here, and throughout this paper, all power series are convergent in a
neighborhood of the origin. The goal of this section is to prove the
following proposition.

\begin{proposition}
\label{pp-1}For any $t$-tuple $(i_{1},\dots ,i_{t})$ of positive integers
with $1\leq i_{1}<\dots <i_{t}\leq s,$ the origin is an isolated fixed point
of $f^{m_{i_{1}}\dots m_{i_{t}}},$%
\begin{equation}
\mu _{f^{m_{i_{1}}\dots m_{i_{t}}}}(0)=(m_{i_{1}}+1)\dots (m_{i_{t}}+1),
\label{index-1}
\end{equation}%
and%
\begin{equation}
P_{m_{i_{1}}\dots m_{i_{t}}}(f,0)=m_{i_{1}}\dots m_{i_{t}}.  \label{index-2}
\end{equation}
\end{proposition}

To prove this result, we first introduce two known results and prove Lemma %
\ref{pp-1-1}.

\begin{lemma}[\protect\cite{LL}]
\label{prod-ind}Let $h_{1}$, $h_{2}\in \mathcal{O}(\Delta ^{n},\mathbb{C}%
^{n}).$ If the origin is an isolated zero of both $h_{1}$ and $h_{2},$ then
the zero order of $h_{1}\circ h_{2}$ at the origin equals to the product of
the zero orders of $h_{1}$ and $h_{2}$ at the origin$,$ say, $\pi
_{h_{1}\circ h_{2}}(0)=\pi _{h_{1}}(0)\pi _{h_{2}}(0).$
\end{lemma}

\begin{theorem}[\textbf{Cronin \protect\cite{Cro}}]
Let $m_{1},\dots ,m_{n}$ be positive integers and let $P=(P_{1},\dots
,P_{n}):\Delta ^{n}\rightarrow \mathbb{C}^{n}$\ be a holomorphic mapping
given by
\begin{equation*}
P_{j}(z_{1},\dots ,z_{n})=\sum_{k=m_{j}}^{\infty }P_{jk}(z_{1},\cdots
,z_{n}),\ j=1,\cdots ,n,
\end{equation*}%
where each $P_{jk}$\ is a homogeneous polynomial of degree $k$\ in $%
z_{1},z_{2},\dots ,z_{n}$. If $0=(0,\dots ,0)$ is an isolated solution of
the system of $n$ equations
\begin{equation*}
P_{jm_{j}}(z_{1},\dots ,z_{n})=0,\ j=1,\dots ,n,
\end{equation*}%
then the origin\ is an isolated zero of the mapping $P$ with order
\begin{equation*}
\pi _{P}(0)=m_{1}m_{2}\dots m_{n}.
\end{equation*}
\end{theorem}

\begin{lemma}
\label{pp-1-1}Let $k$ be a positive integer$.$ Then the $k$-th iteration $%
f^{k}=(f_{1}^{(k)},\dots ,f_{n}^{(k)})$ of $f$ also has the form (\ref%
{form-1}), together with Conditions \ref{c3}---\ref{c5}, say precisely, the
components of $f^{k}$ is given by
\begin{equation}
\left\{
\begin{array}{ll}
f_{j}^{(k)}(x_{1},\dots ,x_{n})= & \lambda _{j}^{k}x_{j}+k\lambda
_{j}^{k-1}x_{j}\sum_{i=1}^{s}a_{ji}x_{i}^{m_{i}}+o_{j}^{(k)},\;1\leq j\leq s,
\\
f_{r}^{(k)}(x_{1},\dots ,x_{n})= & \mu _{r}^{k}x_{r}+R_{r}^{(k)},\;s+1\leq
r\leq n,%
\end{array}%
\right.  \label{pp0}
\end{equation}%
where $\lambda _{1},\dots ,\lambda _{s},\mu _{s+1},\dots ,\mu _{n},a_{ji}$
are the numbers in (\ref{form-1}), each $o_{j}^{(k)}$ satisfies Condition %
\ref{c4} and each $R_{r}^{(k)}$ satisfies Condition \ref{c5}.
\end{lemma}

\begin{proof}
The conclusion is obvious for $r=s+1,\dots ,n.$ By the assumption on the
mapping $f=(f_{1},\dots ,f_{n}),$ for each $j=1,\dots ,$ $s,$ we can write%
\begin{equation*}
f_{j}(x_{1},\dots ,x_{n})=x_{j}\left[ \lambda
_{j}+\sum_{i=1}^{s}a_{ji}x_{i}^{m_{i}}+h_{j}(x_{1}^{m_{1}},\dots
,x_{s}^{m_{s}})\right] +o_{j}^{\prime }(x_{1},\dots ,x_{n}),
\end{equation*}%
where $h_{j}$ and $o_{j}^{\prime }$ are given in Condition \ref{c4}. It is
clear that the terms in the power series $f_{j}^{(k)}(x_{1},\dots ,x_{n})$
affected by $o_{i}^{\prime }\left( x_{1},\dots ,x_{n}\right) $, $i=1,\dots
,s,$ in the iteration process are all of degree $>\left( m_{1}\dots
m_{s}\right) ^{2}.$ Therefore, we may assume that the terms $o_{i}^{\prime
},i=1,\dots ,s,$ are all zero. Then, we can complete the proof by showing
that%
\begin{equation}
f_{j}^{(k)}(x_{1},\dots ,x_{n})=x_{j}\left[ \lambda _{j}^{k}+k\lambda
_{j}^{k-1}\sum_{i=1}^{s}a_{ji}x_{i}^{m_{i}}+h_{j}^{(k)}(x_{1}^{m_{1}},\dots
,x_{s}^{m_{s}})\right] ,  \label{apr5-1}
\end{equation}%
for each $j=1,\dots ,s,$ where $h_{j}^{(k)}(x_{1}^{m_{1}},\dots
,x_{s}^{m_{s}})$ is a polynomial, in $x_{1}^{m_{1}},\dots ,x_{s}^{m_{s}},$
that has the form of (\ref{ad2}). The proof is by induction on $k.$

For $k=1,$ the conclusion is trivial. Assume that (\ref{apr5-1}) holds for $%
k=l$ and $h_{j}^{(l)}$ is in the form of (\ref{ad2}). Then,%
\begin{eqnarray*}
f_{j}^{(l+1)} &=&f_{j}^{(l)}\circ f=f_{j}^{(l)}(f_{1},\dots ,f_{n}) \\
&=&\left[ f_{j}\right] \left[ \lambda _{j}^{l}+l\lambda
_{j}^{l-1}\sum_{i=1}^{s}a_{ji}f_{i}^{m_{i}}+h_{j}^{(l)}((f_{1})^{m_{1}},%
\dots ,(f_{s})^{m_{s}})\right] ,
\end{eqnarray*}%
For each $i=1,\dots ,s,$ considering that $\lambda _{i}^{m_{i}}=1,$ we have
that
\begin{equation*}
(f_{i}(x_{1},\dots ,x_{n}))^{m_{i}}=x_{i}^{m_{i}}+p_{i}(x_{1}^{m_{1}},\dots
,x_{s}^{m_{s}}),
\end{equation*}%
where $p_{i}(x_{1}^{m_{1}},\dots ,x_{s}^{m_{s}})$ is a power series in $%
x_{1}^{m_{1}},\dots ,x_{s}^{m_{s}}$ having the form of (\ref{ad2}) (note
that we assumed $o_{i}^{\prime }\equiv 0$). Thus, it is obvious that there
exist polynomials $q_{j}^{(l+1)}$ and $h_{j}^{(l+1)},$ in $%
x_{1}^{m_{1}},\dots ,x_{s}^{m_{s}},$ having the form of (\ref{ad2}), such
that
\begin{eqnarray*}
f_{j}^{(l+1)}(x_{1},\dots ,x_{n}) &=&\left[ x_{j}\left( \lambda
_{j}+\sum_{i=1}^{s}a_{ji}x_{i}^{m_{i}}+h_{j}(x_{1}^{m_{1}},\dots
,x_{s}^{m_{s}})\right) \right] \\
&&\left[ \lambda _{j}^{l}+l\lambda
_{j}^{l-1}\sum_{i=1}^{s}a_{ji}x_{i}^{m_{i}}+q_{j}^{(l+1)}(x_{1}^{m_{1}},%
\dots ,x_{s}^{m_{s}})\right] \\
&=&x_{j}\left[ \lambda _{j}^{l+1}+\left( l+1\right) \lambda
_{j}^{l}\sum_{i=1}^{s}a_{ji}x_{i}^{m_{i}}+h_{j}^{(l+1)}(x_{1}^{m_{1}},\dots
,x_{s}^{m_{s}})\right] ,
\end{eqnarray*}%
This completes the induction.
\end{proof}

\begin{proof}[\textbf{Proof of Proposition \protect\ref{pp-1}}]
For each $t=1,\dots ,s,$ we first prove (\ref{index-1}) for the $t$-tuple $%
(i_{1},\dots ,i_{t})$ $=$ $(1,$ $2,$ $\dots ,$ $t).$

By Lemma \ref{pp-1-1}, putting $M=m_{1}\dots m_{t}$ and $f^{M}=(f_{1}^{(M)},%
\dots ,f_{n}^{(M)})$, we have, for $j=1,\dots ,s$ and $r=s+1,\dots ,n,$
\begin{equation}
\left\{
\begin{array}{l}
f_{j}^{(M)}(x_{1},\dots ,x_{n})=\lambda _{j}^{M}x_{j}+M\lambda
_{j}^{M-1}x_{j}\sum_{i=1}^{s}a_{ji}x_{i}^{m_{i}}+o_{j}^{(M)}, \\
f_{r}^{(M)}(x_{1},\dots ,x_{n})=\mu _{r}^{M}x_{r}+R_{r}^{(M)},%
\end{array}%
\right.  \label{form-4}
\end{equation}%
where, each $o_{j}^{(M)}=o_{j}^{(M)}(x_{1},\dots ,x_{n})$ satisfies
Condition \ref{c4} and each $R_{r}^{(M)}=R_{r}^{(M)}(x_{1},\dots ,x_{n})$
satisfies Condition \ref{c5}.

By Condition \ref{c1}, we have $\lambda _{j}^{M}=1$\ for $j=1,\dots ,t,$ and
$\lambda _{j}^{M}\neq 1$ for $j=t+1,\dots ,s.$ Then by (\ref{form-4}),
putting $f_{i}^{(M)}=f_{i}^{(M)}(x_{1},\dots ,x_{n})$ for $i=1,\dots ,n,$ we
have%
\begin{equation}
\left\{
\begin{array}{l}
f_{j}^{(M)}-x_{j}=c_{j}x_{j}\sum_{i=1}^{s}a_{ji}x_{i}^{m_{i}}+o_{j}^{(M)},\
\ 1\leq j\leq t, \\
f_{j}^{(M)}-x_{j}=x_{j}(d_{j}+\mathrm{higher\ terms})+o_{j}^{(M)},\ \
t+1\leq j\leq s, \\
f_{r}^{(M)}-x_{r}=e_{r}x_{r}+R_{r}^{(M)},\ \ s+1\leq r\leq n,%
\end{array}%
\right.  \label{big1}
\end{equation}%
where $c_{j}=M\lambda _{j}^{M-1}\neq 0$ for $j=1,\dots ,t,$ $d_{j}=\lambda
_{j}^{M}-1\neq 0$ for $j=t+1,\dots ,s,$ and $e_{r}=\mu _{r}^{M}-1\neq 0$ for
$r=s+1,\dots ,n$ (by Condition \ref{c2})$.$

Let $H:\mathbb{C}^{n}\rightarrow \mathbb{C}^{n}$ be the mapping
\begin{equation*}
(x_{1},\dots ,x_{n})=H(z_{1},\dots ,z_{n})=(z_{1}^{\frac{M}{m_{1}}},\dots
,z_{t}^{\frac{M}{m_{t}}},z_{t+1}^{M},\dots ,z_{s}^{M},z_{s+1},\dots ,z_{n}).
\end{equation*}%
Since for each $j\leq s,o_{j}^{(M)}$ satisfies Condition \ref{c4}$,$ each
term of the power series $o_{j}^{(M)}=o_{j}^{(M)}(x_{1},\dots ,x_{n})$ is
either a monomial of the form
\begin{equation*}
C_{1}x_{j}x_{1}^{m_{1}i_{1}}\dots
x_{t}^{m_{t}i_{t}}x_{t+1}^{m_{t+1}i_{t+1}}\dots x_{s}^{m_{s}i_{s}}\mathrm{\
with\ }i_{1}+\dots +i_{s}\geq 2,
\end{equation*}%
or a monomial of the form
\begin{equation*}
C_{2}x_{1}^{k_{1}}x_{2}^{k_{2}}\dots x_{n}^{k_{n}}\mathrm{\ with\ }%
k_{1}+k_{2}+\dots +k_{n}>(m_{1}\dots m_{s})^{2},
\end{equation*}%
where $C_{1}$ and $C_{2}$ are constants; and therefore, each term of the
power series $o_{j}^{(M)}\circ H(z_{1},\dots ,z_{n})$ is either a monomial
of the form
\begin{equation*}
C_{1}z_{j}^{\alpha _{j}}z_{1}^{Mi_{1}}\dots
z_{t}^{Mi_{t}}z_{t+1}^{Mm_{t+1}i_{t+1}}\dots z_{s}^{Mm_{s}i_{s}}
\end{equation*}%
($\alpha _{j}=\frac{M}{m_{j}}\in \mathbb{N}$ if $1\leq j\leq t,$ or $\alpha
_{j}=M$ if $t+1\leq $ $j\leq s)$ with degree $>2M,$ or a monomial of the
form
\begin{equation*}
C_{2}z_{1}^{\frac{M}{m_{1}}k_{1}}\dots z_{t}^{\frac{M}{m_{t}}%
k_{t}}z_{t+1}^{Mk_{t+1}}\dots z_{s}^{Mk_{s}}z_{s+1}^{k_{s+1}}\dots
z_{n}^{k_{n}}
\end{equation*}%
with degree $>(m_{1}\dots m_{s})^{2}\geq (m_{1}\dots m_{t})^{2}\geq 2M$
(note that all $m_{j}\geq 2$), and then, $o_{j}^{(M)}\circ H(z_{1},\dots
,z_{n})$ is a power series in $z_{1},\dots ,z_{n}$ such that each term has
degree $>2M.$

On the other hand, it is obvious that $R_{r}^{(M)}\circ H(z_{1},\dots
,z_{n}) $ is a power series in $z_{1},\dots ,z_{n}$ consisting of terms of
degree $\geq 2,$ for each $r=s+1,\dots ,n$, by Condition \ref{c5}; and for
each $j=1,\dots ,t,$ it is also obvious that, after applying the
transformation $H, $ $\sum_{i=1}^{s}a_{ji}x_{i}^{m_{i}}$ is changed into
\begin{equation*}
\sum_{i=1}^{t}a_{ji}z_{i}^{M}+\sum_{i=t+1}^{s}a_{ji}z_{i}^{m_{i}M}=%
\sum_{i=1}^{t}a_{ji}z_{i}^{M}+\mathrm{higher\ terms.}
\end{equation*}%
Therefore, putting
\begin{equation}
G_{M}=(g_{1},\dots ,g_{n})=(id-f^{M})\circ H=H-f^{M}\circ H,  \label{ad4}
\end{equation}%
by (\ref{big1}) we have%
\begin{equation}
\left\{
\begin{array}{l}
g_{j}(z_{1},\dots ,z_{n})=-c_{j}z_{j}^{\frac{M}{m_{j}}%
}\sum_{i=1}^{t}a_{ji}z_{i}^{M}+\mathrm{higher\;terms,}\ 1\leq j\leq t, \\
g_{j}(z_{1},\dots ,z_{n})=-d_{j}z_{j}^{M}+\mathrm{higher\;terms,\ }t+1\leq
j\leq s, \\
g_{r}(z_{1},\dots ,z_{n})=-e_{r}z_{r}+\mathrm{higher\;terms,}\ s+1\leq r\leq
n.%
\end{array}%
\right.  \label{big2}
\end{equation}

By Condition \ref{c3}, the $t\times t$ principal submatrix $(a_{ji})_{1\leq
j,i\leq t}$ of the $s\times s$ matrix $(a_{ji})$ is invertible. Thus, $%
0=(0,\dots ,0)$ is an isolated solution of the system of simultaneous
equations:
\begin{equation*}
\left\{
\begin{array}{rrl}
-c_{j}z_{j}^{\frac{M}{m_{j}}}\sum_{i=1}^{t}a_{ji}z_{i}^{M} & =0, & 1\leq
j\leq t, \\
-d_{j}z_{j}^{M} & =0, & t+1\leq j\leq s, \\
-e_{r}z_{r} & =0, & s+1\leq r\leq n.%
\end{array}%
\right.
\end{equation*}

Thus, by Cronin's Theorem and (\ref{big2}), the origin is an isolated zero
of $G_{M}$ with order%
\begin{equation*}
\pi _{G_{M}}(0)=M^{s-t}\prod_{j=1}^{t}\left( \frac{M}{m_{j}}+M\right) =\frac{%
M^{s}}{m_{1}\dots m_{t}}\prod_{j=1}^{t}(1+m_{j})=M^{s-1}%
\prod_{j=1}^{t}(1+m_{j}),
\end{equation*}%
and then the origin is an isolated fixed point $f^{M}.$ On the other hand,
the zero order of $H$ at the origin is $\pi _{H}(0)=\frac{M^{s}}{m_{1}\dots
m_{t}}=M^{s-1}.$ Thus by Lemma \ref{prod-ind} and (\ref{ad4}), we have%
\begin{equation*}
\mu _{f^{m_{1}\dots m_{t}}}(0)=\mu _{f^{M}}(0)=\pi _{id-f^{M}}(0)=\pi
_{G_{M}}(0)/\pi _{H}(0)=\prod_{j=1}^{t}(1+m_{j}).
\end{equation*}%
This implies (\ref{index-1}) for the $t$-tuple $(i_{1},\dots
,i_{t})=(1,\dots ,t).$

Now, let $t$ be any positive integer with $t\leq s.$ We show that (\ref%
{index-1}) holds for any $t$-tuple $(i_{1},\dots ,i_{t})$ of positive
integers with $1\leq i_{1}<\dots <i_{t}\leq s$. Let
\begin{equation*}
(w_{1},\dots ,w_{n})=h(x_{1},\dots ,x_{n})=(x_{\sigma (1)},\dots ,x_{\sigma
(n)})
\end{equation*}%
be a transformation given by a permutation $\sigma $ of $\{1,\dots ,n\}$
with
\begin{equation*}
\sigma (j)=i_{j}\ \mathrm{for\ }j=1,\dots ,t,\ \mathrm{and\ }\sigma (j)=j\
\mathrm{for\ }j=s+1,\dots ,n.
\end{equation*}%
Then in the new coordinate $w=(w_{1},\dots ,w_{n}),$ the original mapping $%
f=(f_{1},\dots ,f_{n})$ is changed into $\widetilde{f}=(\widetilde{f}%
_{1},\dots ,\widetilde{f}_{n})=h\circ f\circ h^{-1}\ $with%
\begin{equation*}
\left\{
\begin{array}{l}
\widetilde{f}_{j}(w_{1},\dots ,w_{n})=\lambda _{\sigma
(j)}w_{j}+w_{j}\sum_{i=1}^{s}a_{\sigma (j)\sigma (i)}w_{i}^{m_{\sigma (i)}}+%
\widetilde{o}_{j},\ j=1,\dots ,s, \\
\widetilde{f}_{r}(w_{1},\dots ,w_{n})=\mu _{r}w_{r}+\widetilde{R}_{r},\
s+1\leq r\leq s,%
\end{array}%
\right.
\end{equation*}%
where $\widetilde{o}_{j}=\widetilde{o}_{j}(w_{1},\dots ,w_{n})=o_{\sigma
(j)}(w_{\sigma ^{-1}(1)},\dots ,w_{\sigma ^{-1}(n)})$ and $\widetilde{R}_{r}=%
\widetilde{R}_{r}(w_{1},\dots ,w_{n})=R_{r}(w_{\sigma ^{-1}(1)},\dots
,w_{\sigma ^{-1}(n)}).$ This is also in the form of (\ref{form-1}) and all
the corresponding Conditions \ref{c1}--\ref{c5} are satisfied, but with a
modification that $\lambda _{1},$ $\dots ,$ $\lambda _{s},$ $m_{1},$ $\dots ,
$ $m_{s},$ $a_{ji},$ $(x_{1},\dots ,x_{s},$ $x_{s+1},\dots x_{n})$ are just
replaced with $\lambda _{\sigma (1)},$ $\dots ,\lambda _{\sigma (s)},$ $%
m_{\sigma (1)},$ $\dots ,$ $m_{\sigma (s)},$ $a_{\sigma (j)\sigma (i)},$ $%
(w_{1},\dots ,w_{s},$ $w_{s+1},\dots ,w_{n}),$ respectively$.$ Then the
above process for the special $t$-tuple $(1,\dots ,t)$ applies to $%
\widetilde{f}$, and then we have
\begin{equation*}
\mu _{\widetilde{f}^{m_{i_{1}}\dots m_{i_{t}}}}(0)=\mu _{\widetilde{f}%
^{m_{\sigma (1)}\dots m_{\sigma (t)}}}(0)=\prod_{j=1}^{t}(m_{\sigma
(j)}+1)=\prod_{j=1}^{t}(m_{i_{j}}+1).
\end{equation*}%
Therefore, by Lemma \ref{cc3.4} we have $\mu _{f^{m_{i_{1}}\dots
m_{i_{t}}}}(0)=\mu _{\widetilde{f}^{m_{i_{1}}\dots
m_{i_{t}}}}(0)=\prod_{j=1}^{t}(m_{i_{j}}+1).$ This completes the proof of (%
\ref{index-1})$.$

To prove (\ref{index-2}), let
\begin{equation*}
\mathfrak{M}_{f}=\{m\in \mathbb{N};\;\mathrm{the\;linear\;part\;of\;}f\;%
\mathrm{at}\;0\;\mathrm{has\;a\ periodic\;point\;of\;period\;}m\}.
\end{equation*}%
Then by Conditions \ref{c1}, \ref{c2}, and Lemma \ref{linear}, $\mathfrak{M}%
_{f}$ is the set consisting of $1$ and all possible products of mutually
distinct numbers taken from the distinct primes $m_{1},\dots ,m_{s},$ and
then by (\ref{index-1}) and Lemma \ref{Th3+1}, for any $t$-tuple $%
(i_{1},\dots ,i_{t})$ of positive integers with $1\leq i_{1}<\dots
<i_{t}\leq s,$ we have
\begin{eqnarray*}
\prod_{j=1}^{t}(1+m_{i_{j}}) &=&\mu _{f^{m_{i_{1}}\dots
m_{i_{t}}}}(0)=\sum_{m\in \mathfrak{M}_{f},m|m_{i_{1}}\dots
m_{i_{t}}}P_{m}(f,0) \\
&=&P_{1}(f,0)+\sum_{k=1}^{t}\sum_{1\leq j_{1}<\dots <j_{k}\leq
t}P_{m_{i_{j_{1}}}\dots m_{i_{j_{k}}}}(f,0),
\end{eqnarray*}%
where the last sum extends over all $k$-tuples $(j_{1},\dots ,j_{k})$ of
positive integers with $1\leq j_{1}<\dots <j_{k}\leq t.$ On the other hand,
by the assumption on $f$, the origin is a simple fixed point of $f$ and so $%
P_{1}(f,0)=\mu _{f}(0)=1.$ Therefore, expanding $%
\prod_{j=1}^{t}(1+m_{i_{j}}) $ and rewriting the above equation, we have%
\begin{equation*}
P_{m_{i_{1}}\dots m_{i_{t}}}(f,0)=m_{i_{1}}\dots
m_{i_{t}}+\sum_{k=1}^{t-1}\sum_{1\leq j_{1}<\dots <j_{k}\leq
t}(m_{i_{j_{1}}}\dots m_{i_{j_{k}}}-P_{m_{i_{j_{1}}}\dots
m_{i_{j_{k}}}}(f,0)).
\end{equation*}

By the above equation, after a standard process of induction on $t\leq s,$
we have $P_{m_{i_{1}}\dots m_{i_{t}}}(f,0)=m_{i_{1}}\dots m_{i_{t}}.$ This
completes the proof.
\end{proof}

\section{\textbf{Normal forms and perturbations}\label{S5}}

In this section we shall combine the normal form method and the perturbation
method to improve Proposition \ref{pp-1}, say, to prove the following
Proposition, in which Condition \ref{c1} is weakened to be (B) and
Conditions \ref{c2}---\ref{c5} are removed.

\begin{proposition}
\label{more}Let $f:\Delta ^{n}\rightarrow \mathbb{C}^{n}$ be a holomorphic
mapping. Assume:

\textrm{(A).} $Df(0)$ is a diagonal matrix.

\textrm{(B). }$Df(0)$ has eigenvalues $\lambda _{1},\dots ,\lambda _{s}$
that are primitive $m_{1}$-th$,\dots ,m_{s}$-th roots of unity,
respectively, and $m_{1},\dots ,m_{s}$ are positive integers.

$\mathrm{(}$\textrm{C}$\mathrm{).}$ The origin is an isolated fixed point of
both $f$ and $f^{M},$ where $M=[m_{1},\cdots ,m_{s}]$ is the least common
multiple of $m_{1},\dots ,m_{s}.$

Then $P_{M}(f,0)>0.$
\end{proposition}

The idea to prove this proposition is to reduce the problem into the easier
one considered in Section \ref{S4}, by small perturbations and coordinate
transformations. We shall perform this in two steps. The first step is to
use the normal form method, Proposition \ref{pp-1} and some results in
Section \ref{S3} to prove Lemma \ref{Key}, which is a weaker version of
Proposition \ref{more}. The second step is to use Lemmas \ref{Key}--\ref%
{number1} and some results in Section \ref{S3} to complete the proof of
Proposition \ref{more}.

The following result is well known in the theory of normal forms (see \cite%
{AP}, p. 84--85).

\begin{lemma}
\label{nor}Let $f:\Delta ^{n}\rightarrow \mathbb{C}^{n}$ be a holomorphic
mapping such that $f(0)=0$ and $Df(0)=(\lambda _{1},\dots ,\lambda _{n})$ is
a diagonal matrix. Then for any positive integer $r,$ there exists a
biholomorphic coordinate transformation in the form of
\begin{equation*}
(y_{1},\dots ,y_{n})=H(x_{1},\dots ,x_{n})=(x_{1},\dots ,x_{n})+\mathrm{%
higher\;terms}
\end{equation*}%
in a neighborhood of the origin such that each component $g_{j}$ of $%
g=(g_{1},\dots ,g_{n})=H^{-1}\circ f\circ H$ has a power series expansion
\begin{equation*}
g_{j}(x_{1},\dots ,x_{n})=\lambda _{j}x_{j}+\sum_{_{i_{1}+\dots
+i_{n}=2}}^{r}c_{i_{1}\dots i_{n}}^{j}x_{1}^{i_{1}}\dots x_{n}^{i_{n}}+%
\mathrm{higher\;terms,\;}j=1,\dots ,n,
\end{equation*}%
in a neighborhood of the origin, where the sum extends over all $n$-tuples $%
(i_{1},\dots ,i_{n})$ of nonnegative integers with
\begin{equation*}
2\leq i_{1}+\dots +i_{n}\leq r\ \mathrm{and\ }\lambda _{j}=\lambda
_{1}^{i_{1}}\dots \lambda _{n}^{i_{n}}.
\end{equation*}
\end{lemma}

\begin{lemma}
\label{eigen}Suppose that $\lambda _{1},\dots ,\lambda _{s}$ satisfy
Condition \ref{c1} and that $(i_{1},\dots ,i_{s})$ is an $s$-tuple of
nonnegative integers. Then for each $j\leq s,$
\begin{equation}
\lambda _{j}=\lambda _{1}^{i_{1}}\dots \lambda _{s}^{i_{s}}  \label{I}
\end{equation}%
if and only if
\begin{equation}
m_{j}|(i_{j}-1),\mathrm{\ }\text{\textrm{and}}\ m_{k}|i_{k}\;\mathrm{%
for\;each\ }k\leq s\;\mathrm{with\;}k\neq j.  \label{II}
\end{equation}
\end{lemma}

\begin{proof}
Recall that Condition \ref{c1} states that $\lambda _{1},\dots ,\lambda _{s}$
are primitive $m_{1}$-th$,\dots ,m_{s}$-th roots of unity, respectively, and
$m_{1},\dots ,m_{s}$ are mutually distinct primes (and thus all $m_{j}\geq 2$%
). Thus, (\ref{II}) implies (\ref{I}) and we can write
\begin{equation*}
\lambda _{k}=e^{\frac{2\pi p_{k}i}{m_{k}}},k=1,\dots ,s,
\end{equation*}%
where $i=\sqrt{-1}$ and each $p_{k}$ is a positive integer with $%
p_{k}<m_{k}. $

For given $j\leq s,$ if (\ref{I}) holds, then putting $n_{j}=i_{j}-1$ and
putting $n_{k}=i_{k}$ for each $k\leq s$ with $k\neq j,$ we have
\begin{equation*}
\lambda _{1}^{n_{1}}\dots \lambda _{s}^{n_{s}}=e^{2\pi i\left( \frac{%
p_{1}n_{1}}{m_{1}}+\dots +\frac{p_{s}n_{s}}{m_{s}}\right) }=1,
\end{equation*}%
and then
\begin{equation*}
l=\frac{p_{1}n_{1}}{m_{1}}+\dots +\frac{p_{s}n_{s}}{m_{s}}
\end{equation*}%
is an integer, and then for each subscript $t$ with $1\leq t\leq s$ we have%
\begin{equation*}
p_{t}n_{t}\prod_{\substack{ j=1  \\ j\neq t}}^{s}m_{j}+m_{t}\sum_{\substack{ %
k=1  \\ k\neq t}}^{s}p_{k}n_{k}\prod_{\substack{ j=1  \\ j\neq k,t}}%
^{s}m_{j}=l\prod_{j=1}^{s}m_{j}.
\end{equation*}%
For each $t\leq s,$ since $p_{t}$ is a positive integer with $p_{t}<m_{t}$
and $m_{1},\dots ,m_{s}$ are distinct primes$,$ $p_{t}\prod_{j=1,j\neq
t}^{s}m_{j}$ and $m_{t}$ are relatively prime. Thus, $m_{t}$ divides $n_{t}$%
. This completes the proof.
\end{proof}

\begin{corollary}
\label{ai}Let $f:\Delta ^{n}\rightarrow \mathbb{C}^{n}$ be a holomorphic
mapping with $f(0)=0$ and assume:

($\mathfrak{A}$). $Df(0)=(\lambda _{1},\dots ,\lambda _{s},\mu _{s+1},\dots
,\mu _{n})$ is a diagonal matrix.

($\mathfrak{B}$). $\lambda _{1},\dots ,\lambda _{s}$ satisfy Condition \ref%
{c1}.

($\mathfrak{C}$). For any $n$-tuple $(i_{1},\dots ,i_{s},i_{s+1},\dots
,i_{n})$ of nonnegative integers,%
\begin{equation*}
\lambda _{j}=\lambda _{1}^{i_{1}}\dots \lambda _{s}^{i_{s}}\mu
_{s+1}^{i_{s+1}}\dots \mu _{n}^{i_{n}}\Longrightarrow i_{s+1}=\dots
=i_{n}=0,\ \mathrm{for\ each\ }j\leq s\mathrm{.}
\end{equation*}

Then, there exists a biholomorphic coordinate transformation in the form of
\begin{equation}
(y_{1},\dots ,y_{n})=H(x_{1},\dots ,x_{n})=(x_{1},\dots ,x_{n})+\mathrm{%
higher\;terms,}  \label{tra}
\end{equation}%
in a neighborhood of the origin, such that the components of $g=(g_{1},\dots
,g_{n})=H^{-1}\circ f\circ H$ are in the form of%
\begin{equation*}
\left\{
\begin{array}{ll}
g_{j}(x_{1},\dots ,x_{n})=\smallskip & \lambda
_{j}x_{j}+x_{j}\sum_{i=1}^{s}b_{ji}x_{i}^{m_{i}}+o_{j},\;j=1,\dots
,s,\smallskip \\
g_{r}(x_{1},\dots ,x_{n})= & \mu _{r}x_{r}+R_{r},\;r=s+1,\dots ,n,%
\end{array}%
\right.
\end{equation*}%
in a neighborhood of the origin, where each $b_{ji}$ is a complex number,
each $o_{j}=o_{j}(x_{1},\dots ,x_{n})$ satisfies Condition \ref{c4} and each
$R_{r}=R_{r}(x_{1},\dots ,x_{n})$ satisfies Condition \ref{c5}$.$
\end{corollary}

\begin{proof}
By ($\mathfrak{A}$) and Lemma \ref{nor}, there exists a biholomorphic
coordinate transformation $H$ in the form of (\ref{tra}) in a neighborhood
of the origin, such that each component of $g=(g_{1},\dots
,g_{n})=H^{-1}\circ f\circ H$ has the expression%
\begin{equation}
g_{j}(x_{1},\dots ,x_{n})=\lambda _{j}x_{j}+\sum_{_{i_{1}+\dots
+i_{n}=2}}^{\left( m_{1}\dots m_{s}\right) ^{2}}c_{i_{1}\dots
i_{s}i_{s+1}\dots i_{n}}^{j}x_{1}^{i_{1}}\dots
x_{s}^{i_{s}}x_{s+1}^{i_{s+1}}\dots x_{n}^{i_{n}}+o_{j}^{\prime }
\label{ad-4.1.2}
\end{equation}%
for\textrm{\ }$j=1,\dots ,s,$ and
\begin{equation*}
g_{r}(x_{1},\dots ,x_{n})=\mu _{r}x_{r}+R_{r},\mathrm{\;for\ }r=s+1,\dots ,n,
\end{equation*}%
in a neighborhood of the origin, where each $o_{j}^{\prime }$ is a power
series in $x_{1},\dots ,x_{n}$ consisting of terms of degree $>\left(
m_{1}\dots m_{s}\right) ^{2}$, each $R_{r}$ is a power series in $%
x_{1},\dots ,x_{n}$ satisfying Condition \ref{c5}, and the sum in (\ref%
{ad-4.1.2}) extends over all $n$-tuples $(i_{1},\dots ,i_{n})$ of
nonnegative integers with%
\begin{equation}
2\leq i_{1}+\dots +i_{n}\leq \left( m_{1}\dots m_{s}\right) ^{2}\ \mathrm{%
and\ }\lambda _{j}=\lambda _{1}^{i_{1}}\dots \lambda _{s}^{i_{s}}\mu
_{s+1}^{i_{s+1}}\dots \mu _{n}^{i_{n}}.  \label{1}
\end{equation}

For each $j\leq s$ and any $n$-tuple $(i_{1},\dots ,i_{n})$ that satisfies (%
\ref{1}), by ($\mathfrak{C}$) we have $i_{s+1}=\dots =i_{n}=0,$ and then by (%
$\mathfrak{B}$) and Lemma \ref{eigen}, we have $m_{j}|(i_{j}-1)$ and $%
m_{k}|i_{k}$ for each $k\leq s$ with $k\neq j.$ Thus, for each $j\leq s,$ (%
\ref{ad-4.1.2}) precisely means that,
\begin{equation*}
g_{j}(x_{1},\dots ,x_{n})=\lambda _{j}x_{j}+x_{j}p_{j}(x_{1}^{m_{1}},\dots
,x_{s}^{m_{s}})+o_{j}^{\prime },
\end{equation*}%
where $p_{j}(x_{1}^{m_{1}},\dots ,x_{s}^{m_{s}})$ is a polynomial in $%
x_{1}^{m_{1}},\dots ,x_{s}^{m_{s}},$ and then we can write%
\begin{equation*}
g_{j}(x_{1},\dots ,x_{n})=\lambda
_{j}x_{j}+x_{j}\sum_{i=1}^{s}b_{ji}x_{i}^{m_{i}}+h_{j}(x_{1}^{m_{1}},\dots
,x_{s}^{m_{s}})+o_{j}^{\prime },
\end{equation*}%
where $h_{j}(x_{1}^{m_{1}},\dots ,x_{s}^{m_{s}})$ is a polynomial, in $%
x_{1}^{m_{1}},\dots ,x_{s}^{m_{s}},$ in the form of (\ref{ad2}). This
completes the proof.
\end{proof}

\begin{lemma}
\label{Key}Let $f:\overline{\Delta ^{n}}\rightarrow \mathbb{C}^{n}$ be a
holomorphic mapping. Assume:

\textrm{(a)} $Df(0)=(\lambda _{1},\dots ,\lambda _{n})$ is a diagonal matrix.

\textrm{(b) }$\lambda _{1},\dots ,\lambda _{s}$ satisfy Condition \ref{c1}: $%
\lambda _{1},\dots ,\lambda _{s}$ are primitive $m_{1}$-th$,\dots ,m_{s}$-th
roots of unity, respectively, and $m_{1},\dots ,m_{s}$ are mutually distinct
primes.

\textrm{(c)} The origin is an isolated fixed point of both $f$ and $%
f^{m_{1}\cdots m_{s}}.$

Then $P_{m_{1}\dots m_{s}}(f,0)>0.$
\end{lemma}

\begin{proof}
It is clear that for the complex numbers $\lambda _{1},\dots ,\lambda _{s},$
the set of all $\left( n-s\right) $-tuples $(\mu _{s+1},\dots ,\mu _{n})$ of
complex numbers so that condition ($\mathfrak{C}$) in Corollary \ref{ai}
fails is a subset of $\mathbb{C}^{n-s}$ with $2(n-s)$-Lebesgue measure zero.
Therefore, for any $\varepsilon >0$, there exist complex numbers $\mu
_{s+1},\dots ,$ $\mu _{n}$ such that%
\begin{equation}
\sum_{k=s+1}^{n}|\mu _{k}-\lambda _{k}|^{2}<\varepsilon ^{2},  \label{b}
\end{equation}%
that Condition \ref{c2} holds for $\mu _{s+1},\dots ,$ $\mu _{n}$, and that
condition ($\mathfrak{C}$) in Corollary \ref{ai} holds for the $n$-tuple $%
(\lambda _{1},\dots ,\lambda _{s},\mu _{s+1},\dots ,\mu _{n})$.

Let $f_{\varepsilon }(z_{1},\dots ,z_{n})$ be the mapping obtained from $%
f(z_{1},\dots ,z_{n})$ by just replacing the linear part of $f$ with $%
(\lambda _{1}z_{1},\dots ,\lambda _{s}z_{s},\mu _{s+1}z_{s+1},\dots ,\mu
_{n}z_{n}),$ say,
\begin{equation*}
f_{\varepsilon }(z_{1},\dots ,z_{n})-f(z_{1},\dots ,z_{n})=(0,\dots
,0,\left( \mu _{s+1}-\lambda _{s+1}\right) z_{s+1},\dots ,\left( \mu
_{n}-\lambda _{n}\right) z_{n}).
\end{equation*}%
Then $f_{\varepsilon }$ satisfies all the assumptions in Corollary \ref{ai}
with
\begin{equation*}
Df_{\varepsilon }(0)=(\lambda _{1},\dots ,\lambda _{s},\mu _{s+1},\dots ,\mu
_{n}),
\end{equation*}%
and then there exists a biholomorphic coordinate transformation in the form
of
\begin{equation*}
(z_{1},\dots ,z_{n})=H_{\varepsilon }(x_{1},\dots ,x_{n})=(x_{1},\dots
,x_{n})+\mathrm{higher\;terms},\ (x_{1},\dots ,x_{n})\in \overline{B},
\end{equation*}%
where $B$ is a ball centered at the origin, such that $g_{\varepsilon
}=(g_{1},\dots ,g_{n})=H_{\varepsilon }^{-1}\circ f_{\varepsilon }\circ
H_{\varepsilon }$ has the expression
\begin{equation}
\left\{
\begin{array}{ll}
g_{j}(x_{1},\dots ,x_{n})=\smallskip & \lambda
_{j}x_{j}+x_{j}\sum_{i=1}^{s}b_{ji}x_{i}^{m_{i}}+o_{j},\;j=1,\dots
,s,\smallskip \\
g_{r}(x_{1},\dots ,x_{n})= & \mu _{r}x_{r}+R_{r},\;r=s+1,\dots ,n,%
\end{array}%
\right.  \label{extr}
\end{equation}%
on $\overline{B}$, where each $b_{ji}$ is a constant, each $o_{j}$ is a
power series satisfying Condition \ref{c4} and each $R_{r}$ is a power
series satisfying Condition \ref{c5}$.$

$g_{\varepsilon }$ is in the form of (\ref{form-1}), together with the
associated Conditions \ref{c1}, \ref{c2}, \ref{c4} and \ref{c5}. Condition %
\ref{c3} may not be satisfied, but for any fixed $g_{\varepsilon },$ by just
modifying the numbers $b_{ji}$ in (\ref{extr}) slightly, we can construct a
sequence $\{g_{k,\varepsilon }\}\subset \mathcal{O}(\overline{B},\mathbb{C}%
^{n}),$ uniformly converging to $g_{\varepsilon }$ on $\overline{B}$ as $%
k\rightarrow \infty ,$ such that each $g_{k,\varepsilon }$ is in the form of
(\ref{form-1}) together with all Conditions \ref{c1}---\ref{c5}.

Now, Proposition \ref{pp-1} applies to each $g_{k,\varepsilon }$, and then
we have%
\begin{equation*}
P_{m_{1}\dots m_{s}}(g_{k,\varepsilon },0)>0,
\end{equation*}%
which implies, by Lemma \ref{cc3.4}, that
\begin{equation*}
P_{m_{1}\dots m_{s}}(H_{\varepsilon }\circ g_{k,\varepsilon }\circ
H_{\varepsilon }^{-1},0)>0.
\end{equation*}%
But it is clear that, for fixed $f_{\varepsilon }$ and $H_{\varepsilon },$ $%
H_{\varepsilon }\circ g_{k,\varepsilon }\circ H_{\varepsilon }^{-1}$
converges to $H_{\varepsilon }\circ g_{\varepsilon }\circ H_{\varepsilon
}^{-1}=f_{\varepsilon },$ uniformly in a neighborhood of the origin as $%
k\rightarrow \infty .$ Thus by Corollary \ref{cc3.5+2}
\begin{equation}
P_{m_{1}\dots m_{s}}(f_{\varepsilon },0)>0,  \label{extr2}
\end{equation}%
provided that the origin is an isolated fixed point of both $f_{\varepsilon
} $ and $f_{\varepsilon }^{m_{1}\dots m_{s}}.$ By (\ref{b}), $f_{\varepsilon
}$ converges to $f$ uniformly in a neighborhood of the origin as $%
\varepsilon \rightarrow 0,$ and then by (c) and Lemma \ref{cc3.1}, for
sufficiently small $\varepsilon ,$ the origin is an isolated fixed point of
both $f_{\varepsilon }$ and $f_{\varepsilon }^{m_{1}\dots m_{s}}$. Hence, (%
\ref{extr2}) holds for sufficiently small $\varepsilon ,$ and hence by (c),
the convergence of $f_{\varepsilon }$ and by Corollary \ref{cc3.5+2},
\begin{equation*}
P_{m_{1}\dots m_{s}}(f,0)>0.
\end{equation*}%
This completes the proof.
\end{proof}

\begin{lemma}
\label{ad-lem4-1}Let $\mathbf{X}$ be a set, let $f:\mathbf{X}\rightarrow
\mathbf{X}$ be a mapping and let $m$ be a positive integer. If $p$ is a
point in $\mathbf{X}$ such that $f^{m}(p)=p,$ then $p$ is a periodic point
of $f$ and the period is a factor of $m.$
\end{lemma}

\begin{proof}
This follows from the definition of periods in Section \ref{S1-1}$.$
\end{proof}

\begin{lemma}
\label{ad-lem4-2}Let $\mathbf{X}$ be a set, let $f:\mathbf{X}\rightarrow
\mathbf{X}$ be a mapping and let $m^{\ast }$ and $n_{1}$ be positive
integers. If $p$ is a point in $\mathbf{X}$ such that
\begin{equation}
f^{m^{\ast }n_{1}}(p)=p,  \label{addl4-2}
\end{equation}%
and $p$ is a periodic point of $f^{m^{\ast }}$ of period $n_{1},$ then $p$
is a periodic point of $f$ with period $mn_{1},$ where $m$ is a factor of $%
m^{\ast }.$
\end{lemma}

\begin{proof}
Assume that $p$ is a periodic point of $f$ of period $L,$ say, $L$ is the
least positive integer such that $f^{L}(p)=p.$ Then, $L|\left( m^{\ast
}n_{1}\right) $ by (\ref{addl4-2}) and Lemma \ref{ad-lem4-1}$,$ and then $L$
can be factorized into $L=mn^{\prime },$ where $m$ and $n^{\prime }$ are
factors of $m^{\ast }$ and $n_{1}$, respectively.

Putting $d=m^{\ast }/m,$ we have
\begin{equation*}
p=f^{L}(p)=f^{dL}(p)=f^{dmn^{\prime }}(p)=f^{m^{\ast }n^{\prime
}}(p)=(f^{m^{\ast }})^{n^{\prime }}(p),
\end{equation*}%
and then $n^{\prime }\geq n_{1},$ for $n_{1}$ is the least positive integer
such that $(f^{m^{\ast }})^{n_{1}}(p)=p.$ Therefore, $n^{\prime }=n_{1}$.
This completes the proof.
\end{proof}

\begin{lemma}
\label{lla}Let $\mathbf{X}$ be a set, let $n_{1},\dots ,n_{s}$ be mutually
distinct primes, let $r_{1},\dots ,r_{s}$ be positive integers, let $f:%
\mathbf{X}\rightarrow \mathbf{X}$ be a mapping and assume that $p\in \mathbf{%
X}$ is a periodic point of $f^{n_{1}^{r_{1}-1}\dots n_{s}^{r_{s}-1}}$ with
period $n_{1}\dots n_{s}$. Then $p$ is a periodic point of $f$ with period $%
n_{1}^{r_{1}}\dots n_{s}^{r_{s}}$.
\end{lemma}

\begin{proof}
By the hypothesis, we have
\begin{equation*}
f^{n_{1}^{r_{1}}\dots n_{s}^{r_{s}}}(p)=\left( f^{n_{1}^{r_{1}-1}\dots
n_{s}^{r_{s}-1}}\right) ^{n_{1}\dots n_{s}}(p)=p.
\end{equation*}%
Thus, $p$ is a periodic point of $f$ and we assume that the period is $L.$
Then by Lemma \ref{ad-lem4-1}, $L$ divides $n_{1}^{r_{1}}\dots
n_{s}^{r_{s}}. $

Since $n_{1},\dots ,n_{s}$ are distinct primes, we have $L=n_{1}^{r_{1}^{%
\prime }}\dots n_{s}^{r_{s}^{\prime }}$ for some nonnegative integers $%
r_{1}^{\prime },\dots ,r_{s}^{\prime }$ with $r_{j}^{\prime }\leq r_{j}$ for
all $j\leq s.$ We must show that $r_{j}^{\prime }=r_{j}$ for all $j\leq s.$

Otherwise, $r_{j}^{\prime }<r_{j}$ for some $j\leq s.$ Without loss of
generality, we assume that $r_{1}^{\prime }<r_{1}.$ Then $%
n_{1}^{r_{1}-1}n_{2}^{r_{2}}\dots n_{s}^{r_{s}}$ is a multiple of $L$ and
then%
\begin{equation*}
f^{n_{1}^{r_{1}-1}n_{2}^{r_{2}}\dots n_{s}^{r_{s}}}(p)=p,
\end{equation*}%
in other words,%
\begin{equation*}
\left( f^{n_{1}^{r_{1}-1}n_{2}^{r_{2}-1}\dots n_{s}^{r_{s}-1}}\right)
^{n_{2}\dots n_{s}}(p)=p,
\end{equation*}%
which contradicts the assumption that $p$ is a periodic point of $%
f^{n_{1}^{r_{1}-1}n_{2}^{r_{2}-1}\dots n_{s}^{r_{s}-1}}$ with period $%
n_{1}n_{2}\dots n_{s}$ (note that all $n_{j}$ are distinct primes and then
all $n_{j}\geq 2).$ This completes the proof.
\end{proof}

\begin{lemma}
\label{number1}Assume that $m_{1},\dots ,m_{s}$ are positive integers such
that for each $j\leq s,$
\begin{equation*}
M=[m_{1},\dots ,m_{s}]>[m_{1},\dots ,m_{j-1},m_{j+1},\dots ,m_{s}].
\end{equation*}%
Then there exist factors $M^{\ast }$ and $M^{\ast \ast }$ of $M$, mutually
distinct primes $n_{1},\dots ,n_{s},$ and positive integers $r_{1},\dots
,r_{s},$ $n_{1}^{\prime },\dots ,n_{s}^{\prime }$, such that for each $j\leq
s,$%
\begin{equation}
n_{j}^{r_{j}}|m_{j}\mathrm{\ but\ }n_{j}^{r_{j}+1}\nmid m_{j},\mathrm{\ and\
}n_{j}^{r_{j}}\nmid m_{k}\mathrm{\ for\ all\ }k\leq s\mathrm{\ with\ }k\neq
j,  \label{00e}
\end{equation}%
\begin{equation}
M=M^{\ast }\prod_{j=1}^{s}n_{j}=M^{\ast \ast }\prod_{j=1}^{s}n_{j}^{r_{j}},
\label{e8}
\end{equation}%
and
\begin{equation}
\frac{M^{\ast }}{m_{j}}=\frac{n_{j}^{\prime }}{n_{j}}\mathrm{\ with\ }%
(n_{j},n_{j}^{\prime })=1,\ j=1,\dots ,s.  \label{ee5}
\end{equation}

Here, $(n_{j},n_{j}^{\prime })$ denotes the greatest common divisor of $%
n_{j} $ and $n_{j}^{\prime }.$
\end{lemma}

\begin{proof}
By the assumption, for each $j\leq s,$ there exists a prime $n_{j}$ and a
positive integer $r_{j}$ such that (\ref{00e}) holds. Then, it is clear that
all the primes $n_{1},\dots ,n_{s}$ are distinct each other. Therefore, $%
n_{1}^{r_{1}}\dots n_{s}^{r_{s}}$ is a factor of $M=[m_{1},\dots ,m_{s}],$
and then $M^{\ast }=[m_{1},\dots ,m_{s}]n_{1}^{-1}\dots n_{s}^{-1}$ and $%
M^{\ast \ast }=[m_{1},\dots ,m_{s}]n_{1}^{-r_{1}}\dots n_{s}^{-r_{s}}$ are
positive integers satisfying (\ref{e8}).

It remains to show the existence of $n_{1}^{\prime },\dots ,n_{s}^{\prime }$
so that (\ref{ee5}) holds. By (\ref{00e}), for each $j\leq s,$ there exists
a positive integer $m_{j}^{\ast \ast }$ such that
\begin{equation}
m_{j}=m_{j}^{\ast \ast }n_{1}^{r_{j1}}\dots n_{s}^{r_{js}}\mathrm{\ with\ }%
r_{jj}=r_{j},  \label{e}
\end{equation}%
and%
\begin{equation}
(m_{j}^{\ast \ast },n_{1}\dots n_{s})=1.  \label{e-0}
\end{equation}%
Then, again by (\ref{00e}), we have for each $j\leq s,$%
\begin{equation}
r_{j}=r_{jj}>r_{kj},\ \mathrm{for\ each\ }k\leq s\ \mathrm{with\ }k\neq j.
\label{adee}
\end{equation}%
Therefore, considering that $n_{1},\dots ,n_{s}$ are mutually distinct
primes, by (\ref{e})--(\ref{adee}) we have%
\begin{equation*}
M=[m_{1},\dots ,m_{s}]=[m_{1}^{\ast \ast },\dots ,m_{s}^{\ast \ast
}]\prod_{j=1}^{s}n_{j}^{r_{jj}}.
\end{equation*}%
This implies that $M^{\ast }=[m_{1}^{\ast \ast },\dots ,m_{s}^{\ast \ast
}]\prod_{j=1}^{s}n_{j}^{r_{jj}-1},$ and then, for each $k\leq s,$ by (\ref{e}%
) we have
\begin{equation*}
\frac{M^{\ast }}{m_{k}}=\frac{[m_{1}^{\ast \ast },\dots ,m_{s}^{\ast \ast
}]\prod_{j=1}^{s}n_{j}^{r_{jj}-1}}{m_{k}^{\ast \ast }n_{1}^{r_{k1}}\dots
n_{s}^{r_{ks}}}=\frac{n_{k}^{\prime }}{n_{k}}
\end{equation*}%
where
\begin{equation*}
n_{k}^{\prime }=\frac{[m_{1}^{\ast \ast },\dots ,m_{s}^{\ast \ast }]}{%
m_{k}^{\ast \ast }}\prod_{\substack{ j=1  \\ j\neq k}}%
^{s}n_{j}^{r_{jj}-r_{kj}-1},
\end{equation*}%
which is a positive integer by (\ref{adee}). Since $n_{1},\dots ,n_{s}$ are
mutually distinct primes, by the previous equality and (\ref{e-0}) we have $%
(n_{k}^{\prime },n_{k})=1.$ This completes the proof.
\end{proof}

\begin{proof}[\textbf{Proof of Proposition \protect\ref{more}}]
Assume that $f:\Delta ^{n}\rightarrow \mathbb{C}^{n}$ is a holomorphic
mapping that satisfies conditions (\textrm{A})--(\textrm{C}) stated in
Proposition \ref{more}:

\textrm{(A).} $Df(0)$ is a diagonal matrix.

\textrm{(B). }$Df(0)$ has eigenvalues $\lambda _{1},\dots ,\lambda _{s}$
that are primitive $m_{1}$-th$,\dots ,m_{s}$-th roots of unity, respectively.

$\mathrm{(}$\textrm{C}$\mathrm{).}$ The origin is an isolated fixed point of
both $f$ and $f^{M},$ where $M=[m_{1},\cdots ,m_{s}].$

We shall show that%
\begin{equation}
P_{M}(f,0)>0.  \label{goal}
\end{equation}%
Since $P_{M}(f,0)$ is completely determined by the local behavior of $f$ at
the origin, we may assume that $f$ is holomorphic on $\overline{\Delta ^{n}}%
, $ by shrinking $\Delta ^{n}.$

By Lemma \ref{cc3.4}, any linear coordinate transformation does not change $%
P_{M}(f,0).$ Therefore, it is without loss of generality to rewrite the
conditions (A) and (B) to be the following conditions (\textrm{D})--(\textrm{%
F}).

(\textrm{D}). $Df(0)=(\lambda _{1},\dots ,\lambda _{s},\lambda _{s+1},\dots
,\lambda _{n})$ is a diagonal matrix.

(\textrm{E}). $\lambda _{1},\dots ,\lambda _{s}$ are primitive $m_{1}$-th$%
,\dots ,m_{s}$-th roots of unity, respectively.

(\textrm{F}). $m_{1},\dots ,m_{s}$ satisfy the assumption in Lemma \ref%
{number1}.

Then all conclusions in Lemma \ref{number1} hold: there exist factors $%
M^{\ast }$ and $M^{\ast \ast }$ of $M$, mutually distinct primes $%
n_{1},\dots ,n_{s},$ positive integers $r_{1},\dots ,r_{s},$ $n_{1}^{\prime
},\dots ,n_{s}^{\prime }$, such that (\ref{00e})--(\ref{ee5}) hold.
Therefore, the following condition also hold.

(\textrm{G}). If $L$ is a positive integer that is the least common multiple
of some numbers in $\{m_{1},\dots ,m_{s}\}$ and is divided by $%
n_{1}^{r_{1}}\dots n_{s}^{r_{s}},$ then $L=M.$

Then, by (\textrm{C}), (\textrm{D})\ and Lemma \ref{cc3.1}, for any positive
number $\varepsilon ,$ by just modifying $\lambda _{s+1},\dots ,\lambda _{n}$
of the linear part of $f$ at the origin$,$ we can construct a holomorphic
mapping $F:\overline{\Delta ^{n}}\rightarrow \mathbb{C}^{n}$ such that the
origin is an isolated fixed point of both $F$ and $F^{M}$,
\begin{equation}
||F-f||_{\overline{\Delta ^{n}}}=\max_{x\in \overline{\Delta ^{n}}%
}|F(x)-f(x)|<\varepsilon ,  \label{sml}
\end{equation}%
and%
\begin{equation}
DF(0)=(\lambda _{1},\dots ,\lambda _{s},\mu _{s+1},\dots ,\mu _{n}),
\label{add1}
\end{equation}%
where $\mu _{s+1},\dots ,\mu _{n}$ are complex numbers with
\begin{equation}
|\mu _{r}|\neq 1,r=s+1,\dots ,n.  \label{add2}
\end{equation}

If we can prove $P_{M}(F,0)>0,$ then by (\ref{sml}), the arbitrariness of $%
\varepsilon $ and Corollary \ref{cc3.5+2}, we shall have (\ref{goal}). In
the rest of the proof, we shall show $P_{M}(F,0)>0.$

Since the origin is an isolated fixed point of $F$ and $F^{M},$ there exists
a ball $B\subset \Delta ^{n}$ centered at the origin with $\overline{B}%
\subset \Delta ^{n},$ such that $F,$ $F^{M^{\ast }}$ and $F^{M}$ are well
defined on $\overline{B},$ that
\begin{equation}
F^{k}(\overline{B})\subset \Delta ^{n}\mathrm{\ for\ all\ }k=1,\dots ,M,
\label{ad-4}
\end{equation}%
and that%
\begin{equation}
\mathrm{Fix}(F|_{\overline{B}})=\mathrm{Fix}(F^{M^{\ast }}|_{\overline{B}})=%
\mathrm{Fix}(F^{M}|_{\overline{B}})=\{0\}.  \label{e0}
\end{equation}

It is clear by (\ref{ee5}) and (\textrm{E}) that $\lambda _{1}^{M^{\ast
}},\dots ,\lambda _{s}^{M^{\ast }}$ are primitive $n_{1}$-th$,\dots ,n_{s}$%
-th roots of unity, respectively, and since $n_{1},\dots ,n_{s}$ are
distinct primes, $\lambda _{1}^{M^{\ast }},\dots ,\lambda _{s}^{M^{\ast }}$
satisfy Condition \ref{c1} for the $s$-tuple $(n_{1},\dots ,n_{s})$. On the
other hand, by (\ref{add1}),%
\begin{equation*}
DF^{M^{\ast }}(0)=(\lambda _{1}^{M^{\ast }},\dots ,\lambda _{s}^{M^{\ast
}},\mu _{s+1}^{M^{\ast }},\dots ,\mu _{n}^{M^{\ast }}).
\end{equation*}%
Therefore, by (\ref{e0}), Lemma \ref{Key} applies to $F^{M^{\ast }}$, and
then%
\begin{equation}
P_{n_{1}\dots n_{s}}(F^{M^{\ast }},0)>0.  \label{r}
\end{equation}

By (\ref{e0}) and Lemma \ref{cc3.3}, there exists a sequence $%
\{F_{j}\}\subset \mathcal{O}(\overline{\Delta ^{n}},\mathbb{C}^{n}),$
converging to $F$ uniformly, such that all fixed points of $F_{j}^{M}$ in $%
\overline{B}$ are hyperbolic.

By (\ref{e0}) and Corollary \ref{cc3.5+2}, if we can prove the following
conclusion (\textrm{H}), we shall obtain the inequality $P_{M}(F,0)>0$
immediately.\medskip

(\textrm{H}). For sufficiently large $j,$ each $F_{j}$ has a periodic point $%
p_{j}$ of period $M$ converging to the origin as $j\rightarrow \infty
.\medskip $

For sufficiently large $j,$ $F_{j}^{M^{\ast }}$ is well defined on $%
\overline{B}$ by (\ref{ad-4}), and uniformly converges to $F^{M^{\ast }}$ on
$\overline{B}$. Thus by (\ref{e0}), (\ref{r}) and Corollary \ref{cc3.2+1}
(iii), for sufficiently large $j,$ $F_{j}^{M}=(F_{j}^{M^{\ast
}})^{n_{1}\dots n_{s}}$ has no fixed point on $\partial B$ and
\begin{equation*}
P_{n_{1}\dots n_{s}}(F_{j}^{M^{\ast }}|_{\overline{B}})=P_{n_{1}\dots
n_{s}}(F^{M^{\ast }},0)>0.
\end{equation*}%
Then, by Lemma \ref{Th3.0} (ii), for sufficiently large $j,$ $F_{j}^{M^{\ast
}}$ has a periodic point $p_{j}\in B$ of period $n_{1}\dots n_{s}.$ And
furthermore, by (\ref{e0}), by the formula
\begin{equation}
F_{j}^{M}(p_{j})=F_{j}^{M^{\ast }n_{1}\dots n_{s}}(p_{j})=\left(
F_{j}^{M^{\ast }}\right) ^{n_{1}\dots n_{s}}(p_{j})=p_{j},  \label{q}
\end{equation}%
and by the fact that $F_{j}^{M}$ uniformly converges to $F^{M}$ on $%
\overline{B},$ we have
\begin{equation}
p_{j}\rightarrow 0\mathrm{\;as\;}j\rightarrow \infty .  \label{qq}
\end{equation}

Thus, by (\ref{e8}), Lemma \ref{lla} and by the formula
\begin{equation*}
\left( F_{j}^{M^{\ast \ast }}\right) ^{n_{1}^{r_{1}-1}\dots
n_{s}^{r_{s}-1}}=F_{j}^{M^{\ast \ast }n_{1}^{r_{1}-1}\dots
n_{s}^{r_{s}-1}}=F_{j}^{M^{\ast }},
\end{equation*}%
$p_{j}$ is a periodic point of $F_{j}^{M^{\ast \ast }}$ with period $%
n_{1}^{r_{1}}\dots n_{s}^{r_{s}}.$ On the other hand, by (\ref{e8}) and (\ref%
{q}) we also have%
\begin{equation*}
F_{j}^{M^{\ast \ast }n_{1}^{r_{1}}\dots
n_{s}^{r_{s}}}(p_{j})=F_{j}^{M}(p_{j})=p_{j}.
\end{equation*}%
Therefore, by Lemma \ref{ad-lem4-2}, $p_{j}$ is a periodic point of $F_{j}$
of period $L_{j}=l_{j}n_{1}^{r_{1}}\dots n_{s}^{r_{s}},$ where $l_{j}\ $is a
positive integer dividing $M^{\ast \ast }.$

Since $l_{j}$ is a bounded sequence of positive integers, by taking
subsequence, we may assume that all $l_{j}$ equal to a fixed positive
integer $l$, say, for each $j,$ $p_{j}$ is a periodic point of $F_{j}$ of
period $L$ with
\begin{equation}
L=ln_{1}^{r_{1}}\dots n_{s}^{r_{s}}.  \label{qq1}
\end{equation}

Now that each $F_{j}$ has a periodic point $p_{j}$ of period $L$ satisfying (%
\ref{qq}) and $F_{j}$ converges to $F$ uniformly$,$ by Corollary \ref{cc3.6}%
, the linear part of $F$ at the origin has a periodic point of period $L,$
and then by (\textrm{E}), (\ref{add1}), (\ref{add2}) and Lemma \ref{linear},
$L$ is the least common multiple of some numbers of $m_{1},\dots ,m_{s}$,
and hence by (\textrm{G}) and (\ref{qq1}), we have
\begin{equation*}
L=M=[m_{1},\dots ,m_{s}],
\end{equation*}%
say, $p_{j}$ is a periodic point of $F_{j}$ of period $M$, and then the
statement (\textrm{H}) is proved. This completes the proof of Proposition %
\ref{more}.
\end{proof}

\section{\textbf{Proof of Theorem \protect\ref{Th1.0}}\label{S6}}

\begin{proof}[\textbf{Proof of Theorem \protect\ref{Th1.0}}]
By shrinking $\Delta ^{n},$ we may assume, without loss of generality, that $%
f$ is holomorphic on $\overline{\Delta ^{n}}$ and the origin is the unique
fixed point of both $f$ and $f^{M}$ in $\overline{\Delta ^{n}}.$ Since the
conclusion in Theorem \ref{Th1.0} trivially holds for\textbf{\ }$M=1,$ we
also assume $M>1.$

If
\begin{equation}
P_{M}(f,0)>0,  \label{pro}
\end{equation}%
then by Corollary \ref{cc3.5+2} (c), there exist a ball $B\subset \Delta ^{n}
$ centered at the origin and a sequence $\{f_{j}\}\subset \mathcal{O}(%
\overline{B},\mathbb{C}^{n})$ such that $f_{j}$ converges to $f|_{\overline{B%
}}$ uniformly and each $f_{j}$ has a periodic point $p_{j}$ of period $M$
with $p_{j}\rightarrow 0$ as $j\rightarrow \infty $. Thus, by Corollary \ref%
{cc3.6} the linear part of $f$ at the origin has a periodic point of period $%
M.$ This proves the necessity in Theorem \ref{Th1.0}.

To prove the sufficiency, assume that the linear part of $f$ at the origin
has a periodic point of period $M.$ Then, by Lemma \ref{linear} there exist
mutually distinct positive integers $m_{1},\dots ,m_{s},\ s\leq n,$ such
that $M=[m_{1},\dots ,m_{s}]$ and $Df(0)$ has eigenvalues $\lambda
_{1},\dots ,\lambda _{s}\ $that are primitive $m_{1}$-th$,\dots ,m_{s}$-th
roots of unity, respectively.

Thus it follows from Corollary \ref{cc3.2+1}, by modifying the linear part
of $f$ slightly, we can construct a sequence of holomorphic mappings $f_{j}:%
\overline{\Delta ^{n}}\rightarrow \mathbb{C}^{n}$ converging to $f$
uniformly, such that for each $j,$ the origin is an isolated fixed point of
both $f_{j}$ and $f_{j}^{M}$, $Df_{j}(0)$ can be diagonalized and $\lambda
_{1},\dots ,\lambda _{s}$ are still eigenvalues of $Df_{j}(0).$ This is
possible because we have assumed that $m_{1},\dots ,m_{s}$ are mutually
distinct, which implies that $\lambda _{1},\dots ,\lambda _{s}$ are
different from each other.

Then, there exists a sequence of invertible linear transformations $H_{j}$
such that
\begin{equation*}
D(H_{j}^{-1}\circ f_{j}\circ H_{j})(0)=(\lambda _{1},\dots ,\lambda _{s},\mu
_{j,s+1},\dots ,\mu _{j,n})
\end{equation*}%
is a diagonal matrix. On the other hand, by Lemma \ref{cc3.4}, the origin is
still an isolated fixed point of both $H_{j}^{-1}\circ f_{j}\circ H_{j}$ and
$(H_{j}^{-1}\circ f_{j}\circ H_{j})^{M}.$ Thus by Proposition \ref{more}, $%
P_{M}(H_{j}^{-1}\circ f_{j}\circ H_{j},0)>0,$ for $\lambda _{1},\dots
,\lambda _{s}$ are primitive $m_{1}$-th, $\dots ,$ $m_{s}$-th roots of unity
and $M=[m_{1},\dots ,m_{s}].$ This implies, by Lemma \ref{cc3.4}, that
\begin{equation*}
P_{M}(f_{j},0)>0.
\end{equation*}%
Hence by the convergence of $f_{j}$ and Corollaries \ref{cc3.5+2} we have $%
P_{M}(f,0)>0.$ This completes the proof of the sufficiency.
\end{proof}

\section{Appendix}

Let $U$ be a bounded open subset of the real vector space $\mathbb{R}^{n}$
containing the origin, $g:\overline{U}\rightarrow \mathbb{R}^{n}$ a
continuous mapping with an isolated fixed point at the origin and let $B$ be
an open ball in $U$ centered at the origin such that $g$ has no other fixed
point in $\overline{B}\backslash \{0\}.$ Then the local\textit{\ Lefschetz
fixed point index} $\mu _{g}(0)$ of $g$ at the origin is defined to be the
topology degree of the mapping
\begin{equation*}
x\mapsto \frac{x-g(x)}{|x-g(x)|}
\end{equation*}%
from the $(n-1)$-sphere $\partial B$ into the unit $(n-1)$-sphere in $%
\mathbb{R}^{n}.$ If, in addition, $g$ is a $C^{1}$ mapping, then%
\begin{equation*}
\mu _{g}(0)=\sum_{x\in B,\ x-g(x)=q}\text{sgn }\det (I-Dg(x)),
\end{equation*}%
where $q\in \mathbb{R}^{n}$ is a regular value (of the mapping $x\longmapsto
x-g(x)$) that is sufficiently close to the origin$,$ $I$ is the $n\times n$
unit matrix and $Dg(x)$ is the differential of $g$ at the origin$,$ which is
identified with the Jacobian matrix (see \cite{LL} and [\ref{SS}])$.$

If the fixed point set \textrm{Fix}$(g)$ of $g$ is a compact subset of $U,$
then one can find a continuous mapping $f:\overline{U}\rightarrow \mathbb{R}%
^{n}$ sufficiently close to $g$ such that \textrm{Fix}$(f)$ is finite and is
contained in $U$, and then one can define the global Lefschetz fixed point
index of $g$ to be%
\begin{equation*}
L(g)=\sum_{p\in U,\ f(p)=p}\mu _{f}(p).
\end{equation*}%
This number is independent of $f,$ when $f$ is close to $g$ enough.

If for a positive integer $M,$ the origin is an isolated fixed point of both
$g$ and $g^{M}$, then the local Dold's index%
\begin{equation*}
P_{M}(g,0)=\sum_{\tau \subset P(M)}(-1)^{\#\tau }\mu _{g^{M:\tau }}(0)
\end{equation*}%
is well defined.

If the fixed point set $\mathrm{Fix}(g^{M})$ of $g^{M}$ is a compact subset
in $U,$ then for each factor $m$ of $M,$ since the fixed point set $\mathrm{%
Fix}(g^{m})$ is a closed subset of $\mathrm{Fix}(g^{M}),$ $\mathrm{Fix}%
(g^{m})$ is again a compact subset of $U$, and then the Lefschetz fixed
point index $L(g^{m})$ is well defined$,$ and so is the Dold's index
\begin{equation*}
P_{M}(g)=\sum_{\tau \subset P(M)}(-1)^{\#\tau }L(g^{M:\tau }).
\end{equation*}%
This number is studied by several authors (see [\ref{Do}], [\ref{FL}], [\ref%
{St}] and [\ref{Za}].) The importance of the number $P_{M}(g)$ is that, when
$P_{M}(g,0)\neq 0,$ any continuous mapping $g_{1}:U\rightarrow \mathbb{R}%
^{n} $ sufficiently close to $g$ has periodic points near the origin with
period $M$ when $M$ is odd, and period $M$ or $M/2$ when $M$ is even,
provided that each fixed point of $g_{1}^{M}$ near the origin is of index $%
+1 $ or $-1$ (see [\ref{FL}])$.$

When $g$ is a holomorphic mapping from $\overline{U}\subset \mathbb{C}^{n}$
into $\mathbb{C}^{n},$ the above definitions for $\mu _{g}(p),P_{M}(g,p)$
and $P_{M}(g)$ agree with the definitions in Sections \ref{S1} and \ref{S3}
(see [\ref{FL}] and \cite{LL}).

\end{document}